\documentclass[psamsfonts]{amsart}

\usepackage{amssymb,amsfonts}
\usepackage[all,arc]{xy}
\usepackage{enumerate}
\usepackage{mathrsfs}
\usepackage{slashbox}
\usepackage{geometry}
\usepackage{lipsum}
\usepackage{tikz}
\usepackage{subfigure}
\usepackage{verbatim}
\usepackage{lipsum}
\usepackage{booktabs} 
\usepackage{xcolor}
\usepackage{subfigure}
\usepackage{graphicx}
\usepackage[msc-links]{amsrefs}
\usepackage{hyperref}
\usepackage{cleveref}
\usepackage{cite}
\usepackage{appendix}
\usepackage{graphics}
\usepackage{graphicx}

\newtheorem{theorem}{Theorem}[section]

\newtheorem{lemma}[theorem]{Lemma}
\newtheorem{proposition}{Proposition}

\theoremstyle{definition}
\newtheorem{definition}[theorem]{Definition}
\newtheorem{remark}{Remark}

\numberwithin{equation}{section}

\title[Hydrodynamic limits on crystal lattice ]{hydrodynamic limits of interacting particle systems on crystal lattices in periodic realizations }

\date{\today}

\begin{document}
\maketitle

\centerline{\scshape Zehao Guan}
\medskip
{\footnotesize
 \centerline{Graduate School of Mathematical Sciences, University of Tokyo, }
   \centerline{Komaba 3-8-1, Meguro-Ku, Tokyo, 153-8914, Japan}

} 


\begin{abstract}
We study the hydrodynamic limits of the simple exclusion processes and the zero range processes on crystal lattices. For a periodic realization of crystal lattice, we derive the hydrodynamic limit for the exclusion processes and the zero range processes, which depends on both the structure of crystal lattice and the periodic realization. Even through the crystal lattices have inhomogeneous local structure, for all periodic realizations, we apply the entropy method to derive the hydrodynamic limits. Also, we discuss how the limit equation depends on the choices of the realizations.
\end{abstract}
\section{Introduction}
The purpose of this paper is to discuss the hydrodynamic limits of interacting particle systems on the crystal lattices. We can regard the interacting particle systems as interacting random walks of lots of particles. Briefly speaking, the hydrodynamic limit, that is to deduce the macroscopic behavior of the system from the microscopic interacting particles, can be regarded as the law of large numbers for these stochastic processes through a proper space-time scaling limit. The limiting macroscopic behavior is described by a deterministic evolution equation, which is called the hydrodynamic equation. The hydrodynamic limits of interacting particle systems have been investigated intensively in the square lattice $\mathbb{Z}^d$, which have their origins in mathematics and physics (See \cite{KL1999}). It is interesting to study the scaling limits of interacting particle systems in much more general spaces. In this direction, Jara \cite{J2009} shows the hydrodynamic limit for the zero range process in the Sierpinski gasket. In \cite{F2008}, Faggionato studies the exclusion process on the percolation clusters.

In this paper, we focus on the crystal lattice, such as the triangular lattice and the hexagonal lattice, which is the simplest extension of the square lattice. The crystal lattice has been studied from the view of discrete geometric analysis by Kotani and Sunada (See \cite{KS2000} \cite{KS2001}). They study random walks on crystal lattices and discuss the relationship between asymptotic behaviors of the random walks and the geometric structures of crystal lattices.
We mention a recent work by Ishiwata, Kawabi and Kotani \cite{IKK2017}, in which they study asymptotic behaviors of non-symmetric random walks on crystal lattices. They establish two kinds of functional central limit theorems for random walks. 

A crystal lattice is defined as an infinite graph $X$ which admits a free action of a free abelian group $\Gamma$ with a finite quotient graph $X_0$. For each positive integer $N$, the subgroup $N\Gamma$ also acts freely on $X$ and we call the finite graph $X_N:=X/N\Gamma$ the $N$-scaling finite graph. Here we consider the exclusion processes and the zero range processes on $X_N$ associated to a symmetric weight function, which can be regarded as a jump rate function. To observe these processes in the continuous space, we embed $X$  into Euclidean space through an embedding map $\Phi$ which respects the group action. We call such an embedding map $\Phi$ a periodic realization of the crystal lattice $X$. For each periodic realization $\Phi$, we construct an embedding map $\Phi_N$ from $X_N$ into a torus such that the image $\Phi_N(X_N)$ converges to a torus as $N$ tends to infinity. (For more details, see Section 2.2) 


Different from the square lattice, the crystal lattice has inhomogeneous local structures, which makes models and the techniques applied much more involved. In \cite{T2012}, Tanaka studied the weakly asymmetric simple exclusion process on the crystal lattices. Tanaka considers the weakly asymmetric simple exclusion process where the weight function is identical to $1$ and a harmonic realization is fixed, then he discusses the influence of the weakly asymmetric part and the macroscopic geometric structures to the macroscopic behaviors of particles. 
We are interested in the case with general symmetric periodic jump rates through periodic realizations (including harmonic realizations).   

In this paper, we deal with the exclusion processes and the zero range processes on crystal lattices and investigate the influence of the geometric structure to the hydrodynamic equation. More precisely, consider the zero range process $\eta(t)$ on $X_N:=(V_N,E_N)$ with generator 
\begin{equation}
L_Nf(\eta):=\sum_{e \in E_N} p(e) g(\eta_{oe}) [f(\eta^e) -f(\eta)],\ \ \ \ \ 
\end{equation}
\noindent where $\eta:=(\eta_x)_{x \in V_N}\in \mathbb{N}^{V_N}$ is the configuration, $g: \mathbb{N} \to \mathbb{R}_+$ with $g(0)=0$ is the jump rate, $\mathbb{N}:=\{0,1,2,\dots,\}$ and $\eta^e$ is given by
\begin{equation*}
\eta^e_x=\left\{
\begin{aligned}
&\eta_{oe}-1\ \ \ \ \ \  x=oe
\\
&\eta_{te}+1\ \ \ \ \ \ x=te
\\
&\eta_x\ \ \ \ \ \ \ \ \ \ \ \  otherwise. 
\end{aligned}
\right.
\end{equation*}
This means, if there are $\eta_{oe}$ particles at site $oe$, independently with the number of particles on other sites, at rate $p(e)g(\eta_{oe})$ one of the particles at $oe$ jumps to $te$. For each periodic realization $\Phi$, define the empirical density by
\begin{equation}
\pi_t^{\Phi,N}(du):=\frac{1}{|V_N|} \sum_{x \in V_N} \eta_x(t) \delta_{\Phi_N(x)} (du).
\end{equation}
\noindent We obtain the behavior of $\pi_t^{\Phi,N}(du)$ as $N \to \infty$ and discuss the influence of the realization $\Phi$ to the hydrodynamic equation. We observe that the diffusion coefficient matrix can be computed by the finite quotient graph $X_0$ and the harmonic realization $\Phi_h$ associated to $\Phi$ (See Theorem \ref{main}). The exclusion process is defined in a similar way (See Section 3.1).  




One of the difficulties in the study of hydrodynamic limit for crystal lattices is that generally there is no gradient type model. According to the types of interactions, interacting particle systems on the square lattice are categorized into the gradient systems and the non-gradient systems. We call the system the gradient system when the current of particles through each bond can be represented by the difference of a local function and its shift. Otherwise, we call the system the non-gradient system.
The gradient condition allows us to do the integration by parts twice to get the hydrodynamic limit equation, which is not possible in the case of non-gradient system (See \cite{KL1999}). For non-gradient systems, Varadhan \cite{V1993} proposed an approach to show the hydrodynamic limits and it has been applied to various non-gradient models. In contrast to the square lattice, for the systems on the crystal lattices, generally we can not perform the integration by parts even once because of the inhomogeneous local structure. The first integration by parts is allowed only when the quotient graph has just one vertex, which is the case for the square lattice. To overcome this difficulty, we utilize the harmonic realization of the crystal lattice. For any periodic realization, there is a unique harmonic realization sharing their macroscopic properties. Hence, to obtain the hydrodynamic limit, we adopt the harmonic realization for which we can replace the time derivative of microscopic particle density by the discrete (weighted) Laplacian of a local function directly without the twice integration by parts.

We look over the outline of the proof for the zero range processes on crystal lattices. Following the entropy method developed in \cite{GPV1988}, we establish the local ergodic theorem, called the replacement lemma. The local ergodic theorem is the key step of the proof since it enables us to replace the local averages by the global averages. The proof of the replacement lemma is based on the one block estimate and the two blocks estimate. 

As an important application of our results, we can deal with a class of non-gradient systems on $\mathbb{Z}^d$ directly without using the non-gradient method. For example, consider the exclusion process on the discrete torus $\mathbb{T}_{ 2N } := \mathbb{ Z } /2N \mathbb{ Z} $, represented by $ \{ 0,1,2,\dots, 2N-1\}$, with generator 

\begin{equation*}
L_Nf(\eta):=\sum_{x  \in \mathbb{ T }_{ 2N }} p(x,x+1) \left\{ f(\eta^{x,x+1}) -f(\eta) \right \},
\end{equation*}



\noindent where

\begin{equation*}
p(x,x+1)=\left\{
\begin{aligned}
&\alpha\ \ \ \ \  x\  \ even
\\
\\
&\beta\ \ \ \ \ x\  \ odd
\end{aligned}
\right.
\end{equation*}
We have that

\begin{equation*}
L_N\eta_x=\left\{
\begin{aligned}
&\beta (\eta_{x-1}-\eta_x)-\alpha (\eta_x-\eta_{x+1}) \ \ \ \ \  x\  \ even
\\
\\
&\alpha (\eta_{x-1}-\eta_x)-\beta (\eta_x-\eta_{x+1})\ \ \ \ \ x\  \ odd
\end{aligned}
\right.
\end{equation*}

\noindent Even though $p(\cdot)$ does not rely on the configuration $\eta$, it turns out to be inhomogeneous and non-gradient. If we consider a new process $ \{ \xi_x:=( \eta_{2x}, \eta_{2x=1} ) \} _{x \in \mathbb{ Z } }$, then $\xi$ is a homogeneous process with more complex state space $ \{ ( 0,0 ), ( 0,1) , ( 1,0 ), ( 1,1 ) \}$. It also might be possible to show the hydrodynamic limit for this model with non-gradient method. 
However, regrading this model as a simple exclusion process on crystal lattice, we can show the hydrodynamic limit directly without using the non-gradient method even though we can not do the integration by parts twice. (See Section 6  Example.1) 


The rest of the paper is organized as follows: In Section 2, we introduce the crystal lattice and construct the $N$-scaling finite graph. In Section 3, we introduce the exclusion process and the zero range process and state our main results. In Section 4, we prove the replacement lemma through the one-block estimate and two-blocks estimate. In Section 5, we discuss how to get the standard realization via the diffusion matrix. In Section 6, we give two examples. In Appendix A, we prove some lemmas. 

\textbf{Notation:} Throughout this paper, $\mathbb{N}=\{0,1,2,\dots,\}$ and $o_N$ means that $o_N \to 0$ as $N \to \infty$.

\section{Crystal Lattice}

\subsection{Crystal lattice}
In this section, we introduce the crystal lattice and fix some notations.

Let $X=(V,E)$ be a locally finite connected graph, where $V$ is the set of vertices and $E$ is the set of all oriented edges. For an oriented edge $e \in E$, we denote by $oe$ the origin of $e$, by $te$ the terminus of $e$ and by $\bar{e}$ the inverse edge of $e$. 
We call $X=(V,E)$ is a $\Gamma$-crystal lattice if a group $\Gamma \cong \mathbb{ Z }^d$ acts on $X$ freely and the quotient graph $X/ \Gamma$ is a finite graph, denoted by $X_0=(V_0,E_0)$. More precisely, each $\sigma \in \Gamma$ defines a graph isomorphism $\sigma : X \to X$ and the graph isomorphism is fixed point-free except for $\sigma =id$. 
Let $p(\cdot)$ be a symmetric $\Gamma$-periodic weight function on $E$, that is, $p(e)=p(\bar{e})>0$ and $p(\sigma e)=p(e)$ for all $e \in E$, $\sigma \in \Gamma$. 
The dimension of $X$, symbolically dim$X$, is defined to be the rank of $\Gamma$. We will embed $X$ into the Euclidean space $\mathbb{R}^d$ of dimension $d=rank \Gamma$. For a fixed $ x_0 \in V$, we can take a fundamental domain $ D_{x_0} \subset V$ such that $ x_0 \in D_{x_0}$ and $D_{x_0}$ is connected in the following sense: For any $x, y \in D_{x_0}$ there exist a path $e_1, . . . , e_l$ in $E$ such that $oe_1 = x,te_l = y$ and $oe_1,te_1,...,oe_l,te_l$ are all in $D_{x_0}$. This kind of set $D_{x_0}$ always exists if we take a spanning tree in $X_0$ and its lift in $X$.

Let $\phi$ be an injective homomorphism: $\Gamma \to \mathbb{R}^d$ such that there exists a basis $u_1,\cdots,u_d \in \mathbb{R}^d$, 
$$\phi(\Gamma)=\{ \sum_{i=1}^d k_i u_i\ | \ k_i\ integers\}.$$
\noindent Together with the vector translation as the group action, we call the image $\phi(\Gamma)$ a lattice group.

\begin{definition}
We call an embedding $\Phi : X \to \mathbb{R}^d$ is a periodic realization if there exists some homomorphism: $\phi:\ \Gamma \to \mathbb{R}^d$ such that $\Phi$ is $\phi$-periodic, i.e., $\Phi(\sigma x)=\Phi(x)+\phi(\sigma)$, for every $x \in V$ and every $\sigma \in \Gamma$.  Furthermore, we call $\Phi$ is a harmonic realization if  $\Phi$ is periodic and harmonic, where ``harmonic" means, $\sum_{e \in E_x}  p(e)[\Phi(te) -\Phi(oe)] =0$, for every $x \in V$, where $E_x:=\{ e \in E | oe=x\}$. 
\end{definition}
In this paper, realizations are always assumed to be periodic. Note that $\phi$ depends on $\Phi$, so we call $\phi(\Gamma)$ the lattice group of $\Phi$. 
Given a periodic realization $\Phi$, define $v(e):=\Phi(te)-\Phi(oe)$ for $e \in E$. By the periodicity, $v$ induces a map on $E_0$, also denoted by $v$.  Define a $d \times d$ symmetric and positive definite matrix by

\begin{equation}
\mathbb{D}_\Phi :=\frac{1}{| V_0 |} \left(\sum_{e \in E_0} p(e) v_i(e)v_j(e) \right)_{i,j=1,\dots,d},
\end{equation}

\noindent which is called the diffusion coefficient matrix of $\Phi$. 

Let $\phi(\Gamma)=\{ \sum_{i=1}^d k_i u_i \}$ be a lattice group, define the fundamental parallelotope for $\phi(\Gamma)$, by setting

\begin{equation}
D_{\phi}:= \left \{ \sum_{i=1}^d t_i u_i\ |\ 0\le t_i<1,i=1,\dots,d \right \}.
\end{equation}
We also define the fundamental parallelotope $D_\Phi$ for realization $\Phi$ by setting $D_\Phi:=D_\phi$, where $\phi(\Gamma)$ is the lattice group of $\Phi$.

 For each periodic realization $\Phi$, define the energy of $\Phi$ by setting
 \begin{equation}
 E(\Phi):=\frac{1}{2} \sum_{e \in E_0} p(e) ||v(e)||^2.
 \end{equation}
 Here $|| \cdot ||$ represents the length of the vector in $\mathbb{ R }^d$.

 \begin{definition}
 We call $\Phi$ is a standard realization if $\Phi$ minimizes the energy among the periodic realizations with fixed volume of fundamental parallelotope, i.e., for all $\Phi^\prime$ with 
 $vol(D_{\Phi^\prime})=vol(D_{\Phi})$, it holds that

 $$E(\Phi) \le E(\Phi^\prime).$$ 
 \end{definition}
 
 \begin{remark} \label{har}
 It has been shown that the harmonic realization minimizes the energy in the family of periodic realizations with the same lattice group and it is unique up to a translation. Furthermore, for a fixed lattice group $\phi ( \Gamma )  $,  the harmonic realization can be obtained by solving the equations 
 
  \begin{equation} \label{harmonic}
  \begin{aligned}
  \sum_{e \in E_x}  p(e)[\Phi(te) -\Phi(oe)] =0, \ \ \  \Phi ( \sigma x )= \Phi ( x ) + \phi ( \sigma ),\ \ \ x \in V_0,\ \ \sigma  \in \phi ( \Gamma ) 
  \end{aligned}
  \end{equation}
 Let $\Phi_{\phi}^h$ be the unique harmonic realization (up to a translation) associated to the lattice group $\phi(\Gamma)$. Thus, to find the standard realization, it suffices to find the lattice group $\phi(\Gamma)$ such that $\Phi_\phi^h$ minimizes the energy with fixed volume. 
  For more details, see \cite{KS2000}, \cite{KS2001}.
  \end{remark}
Let us see some examples of crystal lattices.
\\

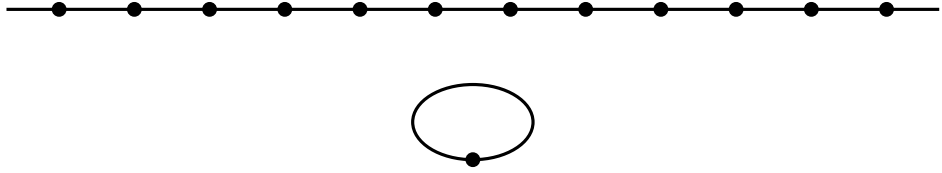
\begin{figure}

      	\begin{tikzpicture}  [xscale = 1, yscale = 1]
	\draw [ black, very thick]  (-3.7, 0) -- (8.7, 0) ; 
	\draw[very thick, fill=black] (-3, 0) circle[radius = 2.2pt] ; 
	\draw[very thick, fill=black] (-2, 0) circle[radius = 2.2pt] ; 
	\draw[very thick, fill=black] (-1, 0) circle[radius = 2.2pt] ; 
	\draw[very thick, fill=black] (0, 0) circle[radius = 2.2pt] ; 
	\draw[very thick, fill=black] (1, 0) circle[radius = 2.2pt] ; 
      	\draw [very thick, fill =black] (2, 0) circle[radius = 2.2pt] ;
	\draw[very thick, fill=black] (3, 0) circle[radius = 2.2pt] ; 
      	\draw [very thick, fill =black] (4, 0) circle[radius = 2.2pt] ;
	\draw[very thick, fill=black] (5, 0) circle[radius = 2.2pt] ; 
      	\draw[very thick, fill=black] (6, 0) circle[radius = 2.2pt] ; 
         \draw[very thick, fill=black] (7, 0) circle[radius = 2.2pt] ; 
         \draw[very thick, fill=black] (8, 0) circle[radius = 2.2pt] ; 
         \draw [black, very thick] (2.5,-1.5) ellipse (0.8 and 0.5); 
         \draw[very thick, fill=black] (2.5,-2) circle[radius=2.2pt];
    \end{tikzpicture}
    \caption{The image of $\Phi$ in Example 1a.}
     \end{figure}

   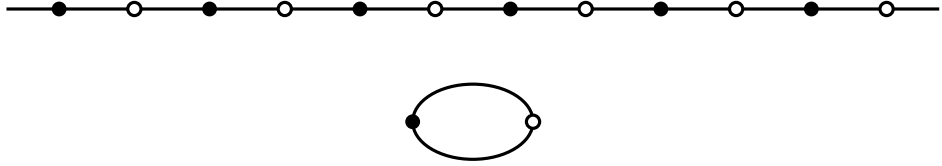
\begin{figure}
      	\begin{tikzpicture} [xscale = 1, yscale = 1]
      	\draw [ black, very thick]  (-3.7, 0) -- (8.7, 0) ; 
	\draw[very thick, fill=black] (-3, 0) circle[radius = 2.2pt] ; 
	\draw[very thick, fill=white] (-2, 0) circle[radius = 2.5pt] ; 
	\draw[very thick, fill=black] (-1, 0) circle[radius = 2.2pt] ; 
	\draw[very thick, fill=white] (0, 0) circle[radius = 2.5pt] ; 
	\draw[very thick, fill=black] (1, 0) circle[radius = 2.2pt] ; 
      	\draw [very thick, fill = white] (2, 0) circle[radius = 2.5pt] ;
	\draw[very thick, fill=black] (3, 0) circle[radius = 2.2pt] ; 
      	\draw [very thick, fill = white] (4, 0) circle[radius = 2.5pt] ;
	\draw[very thick, fill=black] (5, 0) circle[radius = 2.2pt] ; 
      	\draw[very thick, fill=white] (6, 0) circle[radius = 2.5pt] ; 
         \draw[very thick, fill=black] (7, 0) circle[radius = 2.2pt] ; 
         \draw[very thick, fill=white] (8, 0) circle[radius = 2.5pt] ; 
          \draw [black, very thick] (2.5,-1.5) ellipse (0.8 and 0.5); 
        \draw[very thick, fill=black]  (1.7,-1.5) circle[radius=2.2pt];
        \draw[very thick, fill=white]  (3.3,-1.5) circle[radius=2.5pt];
        \end{tikzpicture}
        \caption{The image and the quotient graph of $\Phi$ in Example 1b.}
     \end{figure}

1. One dimensional lattice
\\

1a. The one dimensional standard lattice $X=(V,E)$, where the set of vertices $V=\mathbb{Z}$, the set of edges $E=\{ (x,x+1) , (x+1,x)|\ x \in \mathbb{Z}\}$ 
The group $\mathbb{Z}$ acts freely on $X$ by the additive operation in $\mathbb{Z}$ and the quotient graph consists of one vertex and one loop. We define $\phi : \mathbb{Z} \to \mathbb {R}$ by $\phi(\sigma):=\sigma$ for all $\sigma \in \mathbb{Z}$ and define the embedding map $\Phi(x):=x$ for all $x \in \mathbb{Z}$. Then $\Phi$ is a $\mathbb{Z}$-periodic realization. Furthermore, for the weight function $p(\cdot)$ identically equals to 1, $\Phi$ is harmonic(and also standard) and $\mathbb{D}_\Phi=2$. (See Figure 1)
\\

1b. We give another group action on the above $X$. The group $\mathbb{Z}$ acts freely on $X$ by defining $\sigma x:=x+2\sigma$ for $\sigma \in \mathbb{Z}, x \in V$, then the quotient graph consists of two vertices and two unoriented edges between them. We define $\phi: \mathbb{Z} \to \mathbb{R}$ by $\phi (\sigma ):=2\sigma$ and the embedding $\Phi: X \to \mathbb{R}$ by $\Phi(\sigma 0):=0+\phi(\sigma)$, $\Phi(\sigma 1):=1+\phi(\sigma)$. Then $\Phi$ is a periodic realization. Furthermore, for the weight function $p(\cdot)$ identically equals to 1, $\Phi$ is harmonic(and also standard) and $\mathbb{D}_\Phi=2$. (See Figure 2)
\\

2. The square lattice.
\\

2a. The standard square lattice $X=(V,E)$, where the set of vertices $V=\mathbb{Z}^2$, the set of unoriented edges $E=\{ (x,x+(0,1)),\ (x,x+(1,0)) |\ x \in \mathbb{Z}^2 \}$. Group $\mathbb{Z}^2$ acts freely on $X$ by the additive operation in $\mathbb{Z}^2$ and the quotient graph consists of one vertex and two unoriented loops. We define $\phi : \mathbb{Z}^2 \to \mathbb {R}^2$ by setting $\phi(\sigma):=\sigma$ for $\sigma \in \mathbb{Z}^2$ and define the embedding $\Phi: X \to \mathbb{R}^2$ by 
$\Phi(\sigma (0,0)):= (0,0)+\phi(\sigma)$ for $\sigma \in \mathbb{Z}^2$. Then $\Phi$ is a periodic realization. Furthermore, for the weight function $p(\cdot)$ identically equals to 1, $\Phi$ is harmonic (and also standard) and  
$\mathbb{D}_\Phi=\left(\begin{matrix} 2&0\\ 0&2 \end{matrix}\right)$. (See Figure 3)
\\

2b. We will give another realization for the square lattice $X$. Take a basis $u_1=(1,0), u_2=(1,1)$ in $\mathbb{R}^2$ and define $\phi: \mathbb{Z}^2 \to \mathbb{R}^2$ by $\phi(\sigma):=xu_1+yu_2$ for $\sigma=(x,y) \in \mathbb{Z}^2$. We define the embedding map $\Phi(\sigma (0,0)):=(0,0)+\phi(\sigma)$ for $\sigma \in \mathbb{Z}^2$. Then, for the weight function $p(\cdot)$ identically equals to 1, $\Phi$ is a harmonic realization and $\mathbb{D}_\Phi=\left( \begin{matrix} 4&2\\2&2 \end{matrix}\right)$. (See Figure 3)


\begin{figure}
\begin{tikzpicture} [xscale=0.5,yscale=0.5]
\draw[black, very thick] ( -18.5,-5.5 ) grid ( -8.5, 5.5); 

         \draw[black, very thick] ( -5,-5 ) -- (5 ,5 ); 
         \draw[black, very thick] ( -4,-5 ) -- ( 6, 5); 
         \draw[black, very thick] ( -3,-5 ) -- ( 6, 4); 
         \draw[black, very thick] ( -2,-5 ) -- ( 6, 3); 
         \draw[black, very thick] ( -1,-5 ) -- ( 6, 2); 
         \draw[black, very thick] ( 0,-5 ) -- ( 6, 1); 
         \draw[black, very thick] ( 1,-5 ) -- ( 6, 0); 
         \draw[black, very thick] ( 2,-5 ) -- ( 6, -1); 
         \draw[black, very thick] ( 3,-5) -- ( 6, -2); 
         \draw[black, very thick] ( 4,-5 ) -- ( 6, -3);
         \draw[black, very thick] ( 5,-5 ) -- ( 6, -4);  
         
         \draw[black, very thick] ( -6,-5) -- ( 4, 5); 
         \draw[black, very thick] ( -6,-4 ) -- ( 3, 5); 
         \draw[black, very thick] ( -6,-3 ) -- ( 2, 5); 
         \draw[black, very thick] ( -6,-2 ) -- ( 1, 5); 
         \draw[black, very thick] ( -6,-1 ) -- ( 0, 5); 
         \draw[black, very thick] ( -6,0 ) -- ( -1, 5); 
         \draw[black, very thick] ( -6,1 ) -- ( -2, 5); 
         \draw[black, very thick] ( -6,2 ) -- ( -3, 5); 
         \draw[black, very thick] ( -6,3 ) -- ( -4, 5); 
         \draw[black, very thick] ( -6,4 ) -- ( -5, 5);

          \draw[black, very thick] ( -6,-4.5 ) -- ( 6, -4.5); 
          \draw[black, very thick] ( -6,-3.5 ) -- ( 6, -3.5); 
          \draw[black, very thick] ( -6,-2.5 ) -- ( 6, -2.5); 
          \draw[black, very thick] ( -6,-1.5 ) -- ( 6, -1.5); 
          \draw[black, very thick] ( -6,-0.5 ) -- ( 6, -0.5); 
          \draw[black, very thick] ( -6,0.5 ) -- ( 6, 0.5); 
          \draw[black, very thick] ( -6,1.5 ) -- ( 6, 1.5); 
          \draw[black, very thick] ( -6,2.5 ) -- ( 6, 2.5); 
          \draw[black, very thick] ( -6,3.5 ) -- ( 6, 3.5); 
          \draw[black, very thick] ( -6,4.5 ) -- ( 6, 4.5);

\end{tikzpicture}
\caption{The images of $\Phi$ in Example 2a and 2b. }
\end{figure}
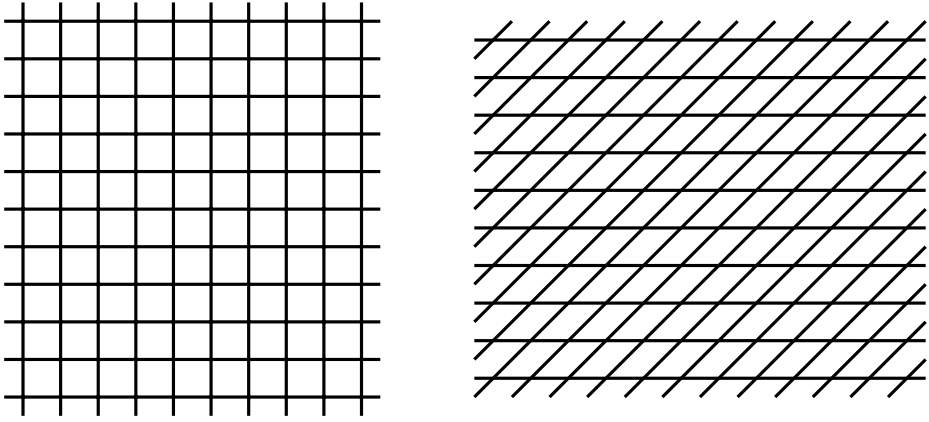

\vspace{0.5cm}

3. The hexagonal lattice. 
\\

3a. The quotient graph of the hexagonal lattice consists of two vertices and three edges as unoriented graph. Define a fundamental graph $D$ by setting the set of vertices $\{x_0,x_1,x_2,x_3\}$ and set of edges $\{e_1,e_2,e_3, \bar {e}_1, \bar {e}_2, \bar { e }_3 \}$, where $e_i:=(x_0,x_i) $ for $i=1,2,3$. The hexagonal lattice is obtained by a $\mathbb{Z}^2$-action on $D$. Take a basis $\{u_1=(\sqrt{3},0), u_2=(\sqrt{3}/2, 3/2)\}$ of $\mathbb{R}^2$ and define $\phi: \mathbb{Z}^2 \to \mathbb{R}^2$ by $\phi(\sigma):=xu_1+yu_2$ for $\sigma=(x,y) \in \mathbb{Z}^2$. We define the embedding $\Phi: X \to \mathbb{R}^2$ by setting $\Phi(\sigma x_0):=(0,0)+\phi(\sigma)$, $\Phi(\sigma x_1):=(0,1)+\phi(\sigma)$, $\Phi(\sigma x_2):=(-\sqrt{3}/2,-1/2)+\phi(\sigma)$, $\Phi(\sigma x_3):=(\sqrt{3}/2,-1/2)+\phi(\sigma)$ for $\sigma \in \mathbb{Z}^2$. Then, for the weight function $p(\cdot)$ identically equals to 1, $\Phi$ is a harmonic realization(and also standard) and $\mathbb{D}_\Phi=\left(\begin{matrix} 3/2& 0\\ 0 &3/2 \end{matrix} \right)$. (See Figure 4)
\\

3b. We consider another realization of the hexagonal lattice. We choose the basis  $\{u_1=(2,0), u_2=(0, 1)\}$ of $\mathbb{R}^2$ and define $\phi: \mathbb{Z}^2 \to \mathbb{R}^2$ by $\phi((x,y)):=xu_1+yu_2$ for $(x,y) \in \mathbb{Z}^2$. We define the embedding $\Phi: X \to \mathbb{R}^2$ by setting $\Phi(\sigma x_0):=(0,0)+\phi(\sigma)$, $\Phi(\sigma x_1):=(0,1)+\phi(\sigma)$, $\Phi(\sigma x_2):=(-1,0)+\phi(\sigma)$, $\Phi(\sigma x_3):=(1,0)+\phi(\sigma)$ for $\sigma \in \mathbb{Z}^2$. Then, for the weight function $p(\cdot)$ identically equals to 1, $\Phi$ is not a harmonic realization. Indeed, for $x=(0,0) \in \mathbb{Z}^2$, $\sum_{e\in E_x} p( e) [\Phi(te)-\Phi(oe)]= (0,1)+ (-1,0)+  (1,0)=( 0,1 )\ne(0,0)$. (See Figure 4)



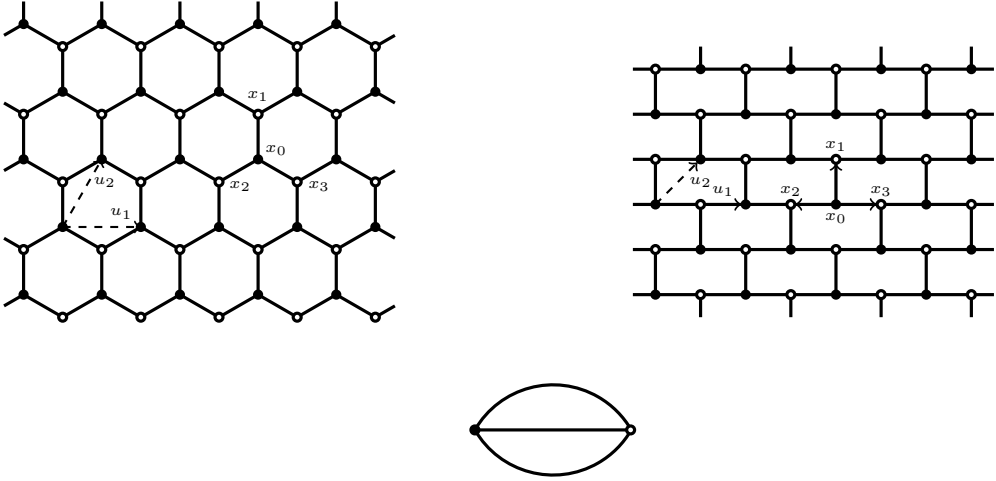
\begin{figure}
     \begin{tikzpicture}[xscale=0.6,yscale=0.6]

      \draw[black, very thick] ( 0,0 ) -- ( {-sqrt(3)/4}, -0.25);   
        \draw[black, very thick] ( 0,0 ) -- ( {sqrt(3)/2}, -0.5); 
        \draw[black, very thick] ( {sqrt(3)},0 ) -- ( {sqrt(3)/2}, -0.5); 
        \draw[black, very thick] ( {sqrt(3)},0 ) -- ( {3*sqrt(3)/2 }, -0.5); 
        \draw[black, very thick] ( {2*sqrt(3)},0 ) -- ( {3*sqrt(3)/2 }, -0.5); 
        \draw[black, very thick] ( {2*sqrt(3)},0 ) -- ( {5*sqrt(3)/2 }, -0.5); 
         \draw[black, very thick] ( {3*sqrt(3)},0 ) -- ( {5*sqrt(3)/2 }, -0.5); 
          \draw[black, very thick] ( {3*sqrt(3)},0 ) -- ( {7*sqrt(3)/2 }, -0.5); 
           \draw[black, very thick] ( {4*sqrt(3)},0 ) -- ( {7*sqrt(3)/2 }, -0.5); 
            \draw[black, very thick] ( {4*sqrt(3)},0 ) -- ( {9*sqrt(3)/2 }, -0.5); 
             \draw[black, very thick] ( {19*sqrt(3)/4},-0.25 ) -- ( {9*sqrt(3)/2 }, -0.5);

              \draw[black, very thick] ( 0,1 ) -- ( {-sqrt(3)/4}, 1.25);   
            \draw[black, very thick] ( 0,1 ) -- ( {sqrt(3)/2}, 1.5); 
             \draw[black, very thick] ( {sqrt(3)},1 ) -- ( {sqrt(3)/2}, 1.5); 
            \draw[black, very thick] ( {sqrt(3)},1 ) -- ( {3*sqrt(3)/2 }, 1.5); 
        \draw[black, very thick] ( {2*sqrt(3)},1 ) -- ( {3*sqrt(3)/2 }, 1.5); 
        \draw[black, very thick] ( {2*sqrt(3)},1 ) -- ( {5*sqrt(3)/2 }, 1.5); 
         \draw[black, very thick] ( {3*sqrt(3)},1 ) -- ( {5*sqrt(3)/2 }, 1.5); 
          \draw[black, very thick] ( {3*sqrt(3)},1 ) -- ( {7*sqrt(3)/2 }, 1.5); 
           \draw[black, very thick] ( {4*sqrt(3)},1 ) -- ( {7*sqrt(3)/2 }, 1.5); 
            \draw[black, very thick] ( {4*sqrt(3)},1 ) -- ( {9*sqrt(3)/2 }, 1.5); 
             \draw[black, very thick] ( {19*sqrt(3)/4},1.25 ) -- ( {9*sqrt(3)/2 }, 1.5); 
                          
              \draw[black, very thick] ( 0,3 ) -- ( {-sqrt(3)/4}, 2.75);   
        \draw[black, very thick] ( 0,3) -- ( {sqrt(3)/2}, 2.5); 
        \draw[black, very thick] ( {sqrt(3)},3 ) -- ( {sqrt(3)/2}, 2.5); 
        \draw[black, very thick] ( {sqrt(3)},3 ) -- ( {3*sqrt(3)/2 }, 2.5); 
        \draw[black, very thick] ( {2*sqrt(3)},3 ) -- ( {3*sqrt(3)/2 }, 2.5); 
        \draw[black, very thick] ( {2*sqrt(3)},3 ) -- ( {5*sqrt(3)/2 }, 2.5); 
         \draw[black, very thick] ( {3*sqrt(3)},3 ) -- ( {5*sqrt(3)/2 }, 2.5); 
          \draw[black, very thick] ( {3*sqrt(3)},3 ) -- ( {7*sqrt(3)/2 }, 2.5); 
           \draw[black, very thick] ( {4*sqrt(3)},3 ) -- ( {7*sqrt(3)/2 }, 2.5); 
            \draw[black, very thick] ( {4*sqrt(3)},3 ) -- ( {9*sqrt(3)/2 }, 2.5); 
             \draw[black, very thick] ( {19*sqrt(3)/4},2.75 ) -- ( {9*sqrt(3)/2 }, 2.5); 
                          
               \draw[black, very thick] ( 0,4 ) -- ( {-sqrt(3)/4}, 4.25);   
            \draw[black, very thick] ( 0,4 ) -- ( {sqrt(3)/2}, 4.5); 
             \draw[black, very thick] ( {sqrt(3)},4 ) -- ( {sqrt(3)/2}, 4.5); 
            \draw[black, very thick] ( {sqrt(3)},4) -- ( {3*sqrt(3)/2 }, 4.5); 
        \draw[black, very thick] ( {2*sqrt(3)},4 ) -- ( {3*sqrt(3)/2 }, 4.5); 
        \draw[black, very thick] ( {2*sqrt(3)},4 ) -- ( {5*sqrt(3)/2 }, 4.5); 
         \draw[black, very thick] ( {3*sqrt(3)},4 ) -- ( {5*sqrt(3)/2 }, 4.5); 
          \draw[black, very thick] ( {3*sqrt(3)},4 ) -- ( {7*sqrt(3)/2 }, 4.5); 
           \draw[black, very thick] ( {4*sqrt(3)},4 ) -- ( {7*sqrt(3)/2 }, 4.5); 
            \draw[black, very thick] ( {4*sqrt(3)},4 ) -- ( {9*sqrt(3)/2 }, 4.5); 
             \draw[black, very thick] ( {19*sqrt(3)/4},4.25 ) -- ( {9*sqrt(3)/2 }, 4.5); 
                           
               \draw[black, very thick] ( 0,6 ) -- ( {-sqrt(3)/4}, 5.75);   
        \draw[black, very thick] ( 0,6) -- ( {sqrt(3)/2}, 5.5); 
        \draw[black, very thick] ( {sqrt(3)},6 ) -- ( {sqrt(3)/2}, 5.5); 
        \draw[black, very thick] ( {sqrt(3)},6 ) -- ( {3*sqrt(3)/2 }, 5.5); 
        \draw[black, very thick] ( {2*sqrt(3)},6 ) -- ( {3*sqrt(3)/2 }, 5.5); 
        \draw[black, very thick] ( {2*sqrt(3)},6 ) -- ( {5*sqrt(3)/2 }, 5.5); 
         \draw[black, very thick] ( {3*sqrt(3)},6 ) -- ( {5*sqrt(3)/2 }, 5.5); 
          \draw[black, very thick] ( {3*sqrt(3)},6 ) -- ( {7*sqrt(3)/2 }, 5.5); 
           \draw[black, very thick] ( {4*sqrt(3)},6 ) -- ( {7*sqrt(3)/2 }, 5.5); 
            \draw[black, very thick] ( {4*sqrt(3)},6 ) -- ( {9*sqrt(3)/2 }, 5.5); 
             \draw[black, very thick] ( {19*sqrt(3)/4},5.75 ) -- ( {9*sqrt(3)/2 }, 5.5); 
                            
              \draw[black, very thick] ( 0,0 ) -- ( 0, 1); 
               \draw[black, very thick] ( 0,3 ) -- ( 0, 4);
                \draw[black, very thick] ( 0,6 ) -- ( 0, 6.5); 
                
                 \draw[black, very thick] ( {sqrt(3)},0 ) -- ( {sqrt(3)}, 1); 
               \draw[black, very thick] ( {sqrt(3)},3 ) -- ( {sqrt(3)}, 4);
                \draw[black, very thick] ( {sqrt(3)},6 ) -- ( {sqrt(3)}, 6.5); 
                
                 \draw[black, very thick] ( {2*sqrt(3)},0 ) -- ( {2*sqrt(3)}, 1); 
               \draw[black, very thick] ( {2*sqrt(3)},3 ) -- ( {2*sqrt(3)}, 4);
                \draw[black, very thick] ( {2*sqrt(3)},6 ) -- ( {2*sqrt(3)}, 6.5); 
                
                 \draw[black, very thick] ( {3*sqrt(3)},0 ) -- ( {3*sqrt(3)}, 1); 
               \draw[black, very thick] ( {3*sqrt(3)},3 ) -- ( {3*sqrt(3)}, 4);
                \draw[black, very thick] ( {3*sqrt(3)},6 ) -- ( {3*sqrt(3)}, 6.5); 
                
                 \draw[black, very thick] ( {4*sqrt(3)},0 ) -- ( {4*sqrt(3)}, 1); 
               \draw[black, very thick] ( {4*sqrt(3)},3 ) -- ( {4*sqrt(3)}, 4);
                \draw[black, very thick] ( {4*sqrt(3)},6 ) -- ( {4*sqrt(3)}, 6.5); 
                
                  \draw[black, very thick] ( {sqrt(3)/2},1.5 ) -- ( {sqrt(3)/2}, 2.5); 
               \draw[black, very thick] ( {sqrt(3)/2},4.5 ) -- ( {sqrt(3)/2}, 5.5);
               
                \draw[black, very thick] ( {3*sqrt(3)/2},1.5 ) -- ( {3*sqrt(3)/2}, 2.5); 
               \draw[black, very thick] ( {3*sqrt(3)/2},4.5 ) -- ( {3*sqrt(3)/2}, 5.5);
                
                 \draw[black, very thick] ( {5*sqrt(3)/2},1.5 ) -- ( {5*sqrt(3)/2}, 2.5); 
               \draw[black, very thick] ( {5*sqrt(3)/2},4.5 ) -- ( {5*sqrt(3)/2}, 5.5);
               
                \draw[black, very thick] ( {7*sqrt(3)/2},1.5 ) -- ( {7*sqrt(3)/2}, 2.5); 
               \draw[black, very thick] ( {7*sqrt(3)/2},4.5 ) -- ( {7*sqrt(3)/2}, 5.5);
               
                \draw[black, very thick] ( {9*sqrt(3)/2},1.5 ) -- ( {9*sqrt(3)/2}, 2.5); 
               \draw[black, very thick] ( {9*sqrt(3)/2},4.5 ) -- ( {9*sqrt(3)/2}, 5.5);
               
                \draw[very thick, fill=black]  (0,0) circle[radius=2.2pt];
                 \draw[very thick, fill=black]  ({sqrt(3)},0) circle[radius=2.2pt];
                   \draw[very thick, fill=black]  ({2*sqrt(3)},0) circle[radius=2.2pt];
                     \draw[very thick, fill=black]  ({3*sqrt(3)},0) circle[radius=2.2pt];
                       \draw[very thick, fill=black]  ({4*sqrt(3)},0) circle[radius=2.2pt];
                       
                         \draw[very thick, fill=white]  ({sqrt(3)/2},-0.5) circle[radius=2.5pt];
                          \draw[very thick, fill=white]  ({3*sqrt(3)/2},-0.5) circle[radius=2.5pt];
                           \draw[very thick, fill=white]  ({5*sqrt(3)/2},-0.5) circle[radius=2.5pt];
                            \draw[very thick, fill=white]  ({7*sqrt(3)/2},-0.5) circle[radius=2.5pt];
                             \draw[very thick, fill=white]  ({9*sqrt(3)/2},-0.5) circle[radius=2.5pt];
                             
                              \draw[very thick, fill=black]  (0,3) circle[radius=2.2pt];
                 \draw[very thick, fill=black]  ({sqrt(3)},3) circle[radius=2.2pt];
                   \draw[very thick, fill=black]  ({2*sqrt(3)},3) circle[radius=2.2pt];
                     \draw[very thick, fill=black]  ({3*sqrt(3)},3) circle[radius=2.2pt];
                       \draw[very thick, fill=black]  ({4*sqrt(3)},3) circle[radius=2.2pt];
                       
                        \draw[very thick, fill=white]  ({sqrt(3)/2},2.5) circle[radius=2.5pt];
                          \draw[very thick, fill=white]  ({3*sqrt(3)/2},2.5) circle[radius=2.5pt];
                           \draw[very thick, fill=white]  ({5*sqrt(3)/2},2.5) circle[radius=2.5pt];
                            \draw[very thick, fill=white]  ({7*sqrt(3)/2},2.5) circle[radius=2.5pt];
                             \draw[very thick, fill=white]  ({9*sqrt(3)/2},2.5) circle[radius=2.5pt];
                             
                               \draw[very thick, fill=black]  (0,6) circle[radius=2.2pt];
                 \draw[very thick, fill=black]  ({sqrt(3)},6) circle[radius=2.2pt];
                   \draw[very thick, fill=black]  ({2*sqrt(3)},6) circle[radius=2.2pt];
                     \draw[very thick, fill=black]  ({3*sqrt(3)},6) circle[radius=2.2pt];
                       \draw[very thick, fill=black]  ({4*sqrt(3)},6) circle[radius=2.2pt];
                       
                        \draw[very thick, fill=white]  ({sqrt(3)/2},5.5) circle[radius=2.5pt];
                          \draw[very thick, fill=white]  ({3*sqrt(3)/2},5.5) circle[radius=2.5pt];
                           \draw[very thick, fill=white]  ({5*sqrt(3)/2},5.5) circle[radius=2.5pt];
                            \draw[very thick, fill=white]  ({7*sqrt(3)/2},5.5) circle[radius=2.5pt];
                             \draw[very thick, fill=white]  ({9*sqrt(3)/2},5.5) circle[radius=2.5pt];

                               \draw[very thick, fill=white]  (0,1) circle[radius=2.5pt];
                                 \draw[very thick, fill=white]  ({sqrt(3)},1) circle[radius=2.5pt];
                                   \draw[very thick, fill=white]  ({2*sqrt(3)},1) circle[radius=2.5pt];
                                     \draw[very thick, fill=white]  ({3*sqrt(3)},1) circle[radius=2.5pt];
                                       \draw[very thick, fill=white]  ({4*sqrt(3)},1) circle[radius=2.5pt];
                                       
                                       \draw[very thick, fill=black]  ({sqrt(3)/2},1.5) circle[radius=2.2pt];
                          \draw[very thick, fill=black]  ({3*sqrt(3)/2},1.5) circle[radius=2.2pt];
                           \draw[very thick, fill=black]  ({5*sqrt(3)/2},1.5) circle[radius=2.2pt];
                            \draw[very thick, fill=black]  ({7*sqrt(3)/2},1.5) circle[radius=2.2pt];
                             \draw[very thick, fill=black]  ({9*sqrt(3)/2},1.5) circle[radius=2.2pt];
                             
                               \draw[very thick, fill=black]  ({sqrt(3)/2},4.5) circle[radius=2.2pt];
                          \draw[very thick, fill=black]  ({3*sqrt(3)/2},4.5) circle[radius=2.2pt];
                           \draw[very thick, fill=black]  ({5*sqrt(3)/2},4.5) circle[radius=2.2pt];
                            \draw[very thick, fill=black]  ({7*sqrt(3)/2},4.5) circle[radius=2.2pt];
                             \draw[very thick, fill=black]  ({9*sqrt(3)/2},4.5) circle[radius=2.2pt];
                             
                               \draw[very thick, fill=white]  (0,4) circle[radius=2.5pt];
                                 \draw[very thick, fill=white]  ({sqrt(3)},4) circle[radius=2.5pt];
                                   \draw[very thick, fill=white]  ({2*sqrt(3)},4) circle[radius=2.5pt];
                                     \draw[very thick, fill=white]  ({3*sqrt(3)},4) circle[radius=2.5pt];
                                       \draw[very thick, fill=white]  ({4*sqrt(3)},4) circle[radius=2.5pt];
                                       
                                       \draw[->, thick, black, dashed] ({sqrt(3)/2},1.5)--({3*sqrt(3)/2},1.5);
                                       \draw [black, very thick] ( {5*sqrt(3)/4}, 1.8) node { \tiny $ u_1 $ } ;
                                         \draw[->, thick, black, dashed] ({sqrt(3)/2},1.5)--({sqrt(3)},3);
                                          \draw [black, very thick] ( 1.8, 2.5) node { \tiny $ u_2 $ } ;
                                          
                                           \draw [black, very thick] ( 5.6, 3.2) node { \tiny $ x_0 $ } ;
                                           
                                            \draw [black, very thick] ( {3*sqrt(3)}, 4.4) node { \tiny $ x_1 $ } ;
                                            
                                             \draw [black, very thick] ( 4.8, 2.4) node { \tiny $ x_2 $ } ;
                                             
                                              \draw [black, very thick] ( 6.55, 2.4) node { \tiny $ x_3 $ } ;

         \draw[black, very thick] ( 13.5,0 ) -- ( 21.5, 0); 
         \draw[black, very thick] ( 13.5,1 ) -- ( 21.5, 1); 
         \draw[black, very thick] ( 13.5,2 ) -- ( 21.5, 2);
         \draw[black, very thick] ( 13.5,3) -- ( 21.5, 3);  
          \draw[black, very thick] ( 13.5,4 ) -- ( 21.5, 4); 
           \draw[black, very thick] ( 13.5,5 ) -- ( 21.5, 5);

             \draw[black, very thick] ( 14,0 ) -- ( 14, 1); 
              \draw[black, very thick] ( 14,2 ) -- ( 14, 3); 
               \draw[black, very thick] ( 14,4 ) -- ( 14, 5); 
               
                \draw[black, very thick] ( 15,-0.5 ) -- ( 15, 0); 
                \draw[black, very thick] ( 15,1 ) -- ( 15, 2);
                \draw[black, very thick] ( 15,3 ) -- ( 15, 4);
                \draw[black, very thick] ( 15,5 ) -- ( 15, 5.5);   
                
                 \draw[black, very thick] ( 16,0 ) -- ( 16, 1); 
              \draw[black, very thick] ( 16,2 ) -- ( 16, 3); 
               \draw[black, very thick] ( 16,4 ) -- ( 16, 5); 
               
                 \draw[black, very thick] ( 17,-0.5 ) -- ( 17, 0); 
                \draw[black, very thick] ( 17,1 ) -- ( 17, 2);
                \draw[black, very thick] ( 17,3 ) -- ( 17, 4);
                \draw[black, very thick] ( 17,5 ) -- ( 17, 5.5);   

                 \draw[black, very thick] ( 18,0 ) -- ( 18, 1); 
              \draw[black, very thick] ( 18,2 ) -- ( 18, 3); 
               \draw[black, very thick] ( 18,4 ) -- ( 18, 5); 
               
                 \draw[black, very thick] (19,-0.5 ) -- ( 19, 0); 
                \draw[black, very thick] ( 19,1 ) -- ( 19, 2);
                \draw[black, very thick] ( 19,3 ) -- ( 19, 4);
                \draw[black, very thick] ( 19,5 ) -- ( 19, 5.5);   
                
                  \draw[black, very thick] ( 20,0 ) -- ( 20, 1); 
              \draw[black, very thick] ( 20,2 ) -- ( 20, 3); 
               \draw[black, very thick] ( 20,4 ) -- ( 20, 5); 
               
                \draw[black, very thick] (21,-0.5 ) -- ( 21, 0); 
                \draw[black, very thick] ( 21,1 ) -- ( 21, 2);
                \draw[black, very thick] ( 21,3 ) -- ( 21, 4);
                \draw[black, very thick] ( 21,5 ) -- ( 21, 5.5);   
                
                  \draw[very thick, fill=black]  (14,0) circle[radius=2.2pt];
                    \draw[very thick, fill=black]  (16,0) circle[radius=2.2pt];
                      \draw[very thick, fill=black]  (18,0) circle[radius=2.2pt];
                        \draw[very thick, fill=black]  (20,0) circle[radius=2.2pt];
                        \draw[very thick, fill=black]  (15,1) circle[radius=2.2pt];
                          \draw[very thick, fill=black]  (17,1) circle[radius=2.2pt];
                            \draw[very thick, fill=black]  (19,1) circle[radius=2.2pt];
                              \draw[very thick, fill=black]  (21,1) circle[radius=2.2pt];
                                \draw[very thick, fill=black]  (14,2) circle[radius=2.2pt];
                                  \draw[very thick, fill=black]  (16,2) circle[radius=2.2pt];
                                    \draw[very thick, fill=black]  (18,2) circle[radius=2.2pt];
                                      \draw[very thick, fill=black]  (20,2) circle[radius=2.2pt];
                                        \draw[very thick, fill=black]  (15,3) circle[radius=2.2pt];
                                          \draw[very thick, fill=black]  (17,3) circle[radius=2.2pt];
                                            \draw[very thick, fill=black]  (19,3) circle[radius=2.2pt];
                                              \draw[very thick, fill=black]  (21,3) circle[radius=2.2pt];
                                                \draw[very thick, fill=black]  (14,4) circle[radius=2.2pt];
                                                 \draw[very thick, fill=black]  (16,4) circle[radius=2.2pt];
                                                   \draw[very thick, fill=black]  (18,4) circle[radius=2.2pt];
                                                   \draw[very thick, fill=black]  (20,4) circle[radius=2.2pt]; 
                                                    \draw[very thick, fill=black]  (15,5) circle[radius=2.2pt];
                                          \draw[very thick, fill=black]  (17,5) circle[radius=2.2pt];
                                            \draw[very thick, fill=black]  (19,5) circle[radius=2.2pt];
                                              \draw[very thick, fill=black]  (21,5) circle[radius=2.2pt];

                                             \draw[very thick, fill=white]  (15,0) circle[radius=2.5pt];
                                               \draw[very thick, fill=white]  (17,0) circle[radius=2.5pt];
                                                 \draw[very thick, fill=white]  (19,0) circle[radius=2.5pt];
                                                   \draw[very thick, fill=white]  (21,0) circle[radius=2.5pt];
                                                    \draw[very thick, fill=white]  (15,2) circle[radius=2.5pt];
                                               \draw[very thick, fill=white]  (17,2) circle[radius=2.5pt];
                                                 \draw[very thick, fill=white]  (19,2) circle[radius=2.5pt];
                                                   \draw[very thick, fill=white]  (21,2) circle[radius=2.5pt];
                                                    \draw[very thick, fill=white]  (15,4) circle[radius=2.5pt];
                                               \draw[very thick, fill=white]  (17,4) circle[radius=2.5pt];
                                                 \draw[very thick, fill=white]  (19,4) circle[radius=2.5pt];
                                                   \draw[very thick, fill=white]  (21,4) circle[radius=2.5pt];
                                                    \draw[very thick, fill=white]  (14,1) circle[radius=2.5pt];
                                               \draw[very thick, fill=white]  (16,1) circle[radius=2.5pt];
                                                 \draw[very thick, fill=white]  (18,1) circle[radius=2.5pt];
                                                   \draw[very thick, fill=white]  (20,1) circle[radius=2.5pt];
                                                    \draw[very thick, fill=white]  (14,3) circle[radius=2.5pt];
                                               \draw[very thick, fill=white]  (16,3) circle[radius=2.5pt];
                                                 \draw[very thick, fill=white]  (18,3) circle[radius=2.5pt];
                                                   \draw[very thick, fill=white]  (20,3) circle[radius=2.5pt];
                                                    \draw[very thick, fill=white]  (14,5) circle[radius=2.5pt];
                                               \draw[very thick, fill=white]  (16,5) circle[radius=2.5pt];
                                                 \draw[very thick, fill=white]  (18,5) circle[radius=2.5pt];
                                                   \draw[very thick, fill=white]  (20,5) circle[radius=2.5pt];
                                                   
                                                   \draw [  very thick, black] (10 , -3) arc [start angle = 150, end angle =30 , radius = 2];
                                                     \draw [  very thick, black] (10 , -3) arc [start angle = -150, end angle =-30 , radius = 2];
                                                     \draw[ very thick, black] (10,-3)--(13.4,-3);
                                                      \draw[very thick, fill=black]  (10,-3) circle[radius=2.5pt];
                                                 
                                                       \draw[very thick, fill=white]  (13.45,-3) circle[radius=2.5pt];

                                                         \draw [black, very thick] ( 18, 1.7) node { \tiny $ x_0 $ } ;
                                                           \draw [black, very thick] ( 19, 2.3) node { \tiny $ x_3 $ } ;
                                                             \draw [black, very thick] ( 17, 2.3) node { \tiny $ x_2 $ } ;
                                                               \draw [black, very thick] ( 18, 3.3) node { \tiny $ x_1 $ } ;
                                                               \draw[->,thick, black] (18,2)--(17.1,2);
                                                                 \draw[->,thick, black] (18,2)--(18.9,2);
                                                                   \draw[->,thick, black] (18,2)--(18,2.9);
                                                                     \draw[->,thick, black,dashed] (14,2)--(15.9,2);
                                                                       \draw[->,thick, black,dashed] (14,2)--(14.93,2.93);
                                                                        \draw [black, very thick] ( 15.5, 2.3) node { \tiny $ u_1 $ } ;
                                                                         \draw [black, very thick] ( 15, 2.5) node { \tiny $ u_2 $ } ;

     \end{tikzpicture}
\caption{The images and quotient graph of $\Phi$ in Example 3a an 3b.}
\end{figure}


\subsection{N-scaling finite graph}

Recall that $\Gamma$ is isomorphic to $\mathbb{Z}^d$. For every positive integer $N \ge 1$, $N\Gamma$ is isomorphic to $N \mathbb{Z}^d$. The subgroup $N\Gamma$ acts also freely on $X$ and its quotient graph $X/ N\Gamma$ is also a finite graph, denoted by $X_N=(V_N,E_N)$. Then $\Gamma_N:= \Gamma / N\Gamma \cong \mathbb{Z}^d / N \mathbb{Z}^d$ acts freely on $X_N$. We call $X_N$ the $N$-scaling finite graph. Since $\Phi$ is periodic, the map

$$\frac{1}{N}\Phi: X \to \mathbb{R}^d,$$
satisties $\frac{1}{N} \Phi(\sigma^N x)=\frac{1}{N} \Phi(x) +\phi(\sigma)$, where $\phi(\Gamma)$ is the lattice group of $\Phi$. Let $\mathbb{T}_\phi^d:=\mathbb{R}^d /\phi(\Gamma)$, equipped with the flat metric induced from the Euclidean metric. Then the map $\frac{1}{N} \Phi$ induces the map

$$\Phi_N: X_N \to \mathbb{T}_\phi^d.$$
We call $\Phi_N$ the $N$-scaling map. We can think about $\Phi_N$ as a discrete approximation of the continuous torus $\mathbb{T}^d_\phi$.


\section{Interacting particle systems on crystal lattices}

\subsection{Simple exclusion process on crystal lattices}
Let $X_N:=(V_N,E_N)$ be the $N$-scaling finite graph. Let $Z_N= \{ 0, 1 \} ^ {V_N}$ be the configuration space and $\eta=\{\eta_x\}_{x\in V_N} \in Z_N$ be the configuration. Let $\nu_\rho^N(0\le \rho \le1)$ be the product Bernoulli measure with density $\rho$, i.e., $\nu_\rho^N ( \eta_x=1)= \rho.$ The generator acting on $L^2(Z_N, \nu_\rho^N)$ as 

\begin{equation}
L_Nf(\eta):=\sum_{e \in E_N} p(e) [f(\eta^e) -f(\eta)],\ \ \ \ \ \ \ f \in L^2(Z_N, \nu_\rho^N),
\end{equation}
where 
\begin{equation*}
\eta^e_x=\left\{
\begin{aligned}
&\eta_{te}\ \ \ \ \ \  x=oe
\\
&\eta_{oe}\ \ \ \ \ \ x=te
\\
&\eta_x\ \ \ \ \ \  \  otherwise,
\end{aligned}
\right.
\end{equation*}

\noindent defines a Markov process $\eta(t)$ on $Z_N$ called the simple exclusion process. 

For an arbitrary fixed time $T>0$, let $D([0,T],Z_N) $ be the path space and for a probability measure $\mu^N$ on $Z_N$, let $P_{\mu^N}$ be the distribution on $D([0,T],Z_N) $ of the continuous Markov process with generator $N^2L_N$ and initial measure $\mu^N$ and $E_{\mu^N}$ be the expectation with respect to $P_{\mu^N}$. For the exclusion process, we have the following result:

\begin{theorem} \label{exclusion}
Let $\Phi$ be a periodic realization with lattice group $\phi(\Gamma)$ and $\Phi_h:=\Phi_\phi^h$ be the harmonic realization associated to the lattice group $\phi(\Gamma)$. Let $\rho_0: \mathbb{T}^d_\phi \to [0,1] $ be a measurable function. Assume that the initial measures $\{\mu^N\}$ satisfy that 

\begin{equation}
\varlimsup_{N \to \infty} \mu^N\left[ \left|\frac{1}{|V_N|} \sum_{x \in V_N} G(\Phi_N(x)) \eta_x -\int_{\mathbb{T}_\phi^d} G(u) \rho_0(u) \frac{du}{vol({\mathbb{T}^d_\phi})}\right|>\delta\right]=0
\end{equation}
for every $\delta>0$ and every continuous function $G: \mathbb{T}^d_\phi \to \mathbb{R}$, then for every $t>0$,   
\begin{equation}
\varlimsup_{N \to \infty} P_{\mu^N } \left[ \left|\frac{1}{|V_N|} \sum_{x \in V_N} G(\Phi_N(x)) \eta_x(t) -\int_{\mathbb{T}_\phi^d} G(u) \rho(t,u) \frac{du}{vol(\mathbb{T}^d_\phi)}\right|>\delta\right]=0
\end{equation}
for every $\delta>0$ and every continuous function $G: \mathbb{T}^d_\phi \to \mathbb{R}$, where $\rho(t,u)$ is the unique weak solution of the following linear heat equation

 \begin{equation} \label{eq1}
 \left\{
 \begin{aligned}
  &\frac{\partial }{\partial t} \rho =\nabla \mathbb{D}_{\Phi_h} \nabla  \rho,
  \\
 & \rho(0,\cdot)=\rho_0(\cdot).
  \end{aligned}
  \right.
  \end{equation}
\end{theorem}

In \cite{T2012}, Tanaka obtains the hydrodynamic limit for the simple exclusion process in the case where the weight function $p ( \cdot )$ is identical to 1 and the realization $\Phi $ is harmonic associated with $p(\cdot) \equiv1$. In Theorem \ref{exclusion}, we obtain the hydrodynamic limit for the exclusion process in the case of general periodic realizations and general symmetric periodic weight functions $p( \cdot)$. The proof has two parts. Firstly, we obtain the hydrodynamic limit when $ \Phi $ is harmonic associated with a given symmetric periodic weight function $p( \cdot )$. The proof of this part is similar to \cite{ T2012}, and also similar to the zero-range case we discuss later in the paper, so we omit it here. Secondly, for a general symmetric periodic weight function $p( \cdot ) $, we extend the result to the case of general periodic realization from harmonic realization. This part is proved by the following proposition.


\begin{proposition} \label{be}
Let $ \Phi $ be a periodic realization with lattice group $\phi(\Gamma)=\{ \sum_{i=1}^d k_i u_i\}$.  
 Assume for every $t>0$, it holds that, 

\begin{equation}
\varlimsup_{N \to \infty} P_{\mu^N } \left[ \left|\frac{1}{|V_N|} \sum_{x \in V_N} G(\Phi_N(x)) \eta_x(t) -\int_{\mathbb{T}_\phi^d} G(u) \rho(t,u) \frac{du}{vol(\mathbb{T}^d_\phi)}\right|>\delta\right]=0
\end{equation}
for every $\delta>0$ and every continuous function $G: \mathbb{T}^d_\phi \to \mathbb{R}$, where $\rho(t,u)$ is the unique weak solution of the following linear heat equation

 \begin{equation} \label{eq2}
 \left\{
 \begin{aligned}
  &\frac{\partial }{\partial t} \rho =\nabla \mathbb{D} \nabla  \rho,
  \\
 & \rho(0,\cdot)=\rho_0(\cdot)
  \end{aligned}
  \right.
  \end{equation}
where $\mathbb{D}$ is the diffusion matrix of $ \Phi $. Then, for any periodic realization $ \widetilde { \Phi } $ with lattice group $\widetilde{\phi}(\Gamma)=\{\sum_{i=1}^d k_i \widetilde{u}_i\}$ and every $t>0$, we have

\begin{equation}
\varlimsup_{N \to \infty} P_{\mu^N } \left[ \left|\frac{1}{|V_N|} \sum_{x \in V_N} G(\widetilde{ \Phi }_N(x)) \eta_x(t) -\int_{\mathbb{T}_{\widetilde{\phi}}^d} G(u) \widetilde{ \rho } (t,u) \frac{du}{vol(\mathbb{T}^d_{ \widetilde{ \phi }})}\right|>\delta\right]=0
\end{equation}
for every $\delta>0$ and every continuous function $G: \mathbb{T}^d_{ \widetilde{ \phi }} \to \mathbb{R}$, where $ \widetilde{ \rho }(t,u)$ is the unique weak solution of the following linear heat equation

 \begin{equation} \label{eq3}
 \left\{
 \begin{aligned}
  &\frac{\partial }{\partial t}  \widetilde{ \rho } =\nabla A \mathbb{D} A^T \nabla \widetilde{ \rho },
  \\
 & \widetilde {\rho}(0,\cdot)=\rho_0(A^{-1}\cdot).
  \end{aligned}
  \right.
  \end{equation}
  Here $A$ is the basis transformation matrix from $\{u_1,\dots,u_d\}$ to $\{\widetilde{u}_1,\dots,\widetilde{u}_d\} $.

\end{proposition}

\noindent Proof: First, we note that a realization $\Phi$ is uniquely determined by its values on $D_{x_0}$ and its lattice group, i.e., $\{\Phi (x) , x \in D_{x_0}\}$ and $\phi(\Gamma)=\{ \sum_{i=1}^d k_i u_i\}$. For every $G \in C^2(\mathbb{T}_{ \widetilde { \phi} }^d)$, 


\begin{align*}
\frac{1}{ | V_N | } \sum_{x \in V_N} G(\widetilde{\Phi}_N(x)) \eta_x &= \frac{1}{|V_N|} \sum_{x \in D_{x_0}} \sum_{\sigma \in \Gamma_N} G\left ( \frac{A \phi(\sigma) +\widetilde{\Phi} (x)} {N} \right) 
\eta_{\sigma x}
\\
&=\frac{1}{|V_N|} \sum_{x \in D_{x_0}} \sum_{\sigma \in \Gamma_N} \widetilde{G} \left (\frac{ \phi(\sigma) +A^{-1}\widetilde{\Phi} (x)} {N} \right ) \eta_{\sigma x},
\end{align*}
 where $\widetilde{G}(x):=G(Ax)$ is a function on $\mathbb{T}_{ \phi }^d$. Notice that

 \begin{align*}
&\left| \frac{1}{| V_N |} \sum_{x \in D_{x_0}} \sum_{\sigma \in \Gamma_N} \widetilde{G}  \left (\frac{ \phi(\sigma) +A^{-1}\widetilde{\Phi} (x)} {N} \right ) \eta_{\sigma x} -\frac{1}{|V_N|} \sum_{x \in D_{x_0}} \sum_{\sigma \in \Gamma_N} \widetilde{G} \left (\frac{ \phi(\sigma) +\Phi (x)} {N} \right)\eta_{\sigma x}\right| 
\\
&\ \ \ \ \ \ \ \ \ \ \ \ \le \sup_{x \in D_{x_0}\atop  \sigma \in \Gamma_N} \left|\widetilde{G} \left (\frac{ \phi(\sigma) +A^{-1}\widetilde{\Phi} (x)} {N} \right )-\widetilde{G} \left (\frac{ \phi(\sigma) +\Phi (x)} {N} \right)\right| \frac{1}{|V_N|} \sum_{x \in V_N} \eta_x,
 \end{align*}

\noindent which converges to $0$ in probability as $N \to \infty$. By assumption, we have that

\begin{align*}
\frac{1}{|V_N|} \sum_{x \in V_N} \widetilde{G} (\Phi_N(x)) \eta_x(t) &\to \int_{\mathbb{T}^d_\phi} \widetilde{G}(u) \rho_0(u)  \frac{du}{vol(\mathbb{T}^d_\phi)}
\\
&\quad +\int_0^t ds \int_{\mathbb{T}^d_\phi} \nabla \mathbb{D} \nabla \widetilde{G}(u) \rho(s,u)  \frac{du}{vol(\mathbb{T}^d_\phi)}
\\
&=\int _{ \mathbb{T}^d_{\widetilde{\phi}}} G(u) \rho_0(A^{-1}u)  \frac{1}{|A|}  \frac{du}{vol(\mathbb{T}^d_\phi)} \\
& \quad +\int_0^t \int_{\mathbb{T}^d_{\widetilde{\phi}}} \nabla A \mathbb{D} A^T \nabla G(u) \rho(s, A^{-1} u) \frac{1}{|A|}  \frac{du}{vol(\mathbb{T}^d_\phi)} 
\\
&=\int_{\mathbb{T}^d_{\widetilde{\phi}}} G(u) \widetilde{\rho}_0(u)  \frac{du}{vol(\mathbb{T}^d_{\widetilde{\phi}})} 
\\
&\quad+\int_0^t ds \int_{\mathbb{T}^d_{\widetilde{\phi}}} \nabla A \mathbb{D} A^T\nabla G(u) \widetilde{\rho}(s,u) \frac{du}{vol(\mathbb{T}^d_{\widetilde{\phi}})} 
\end{align*}
in probability as $N \to \infty$. By the triangular inequality, the proof is completed.

\vspace{0.2cm}

Proof of Theorem \ref{exclusion}: For any given periodic realization $ \Phi $, we take the harmonic realization $ \Phi_h $ whose lattice group coincides with $ \Phi $. Firstly, we obtain the hydrodynamic limit of the exclusion process for $ \Phi_h$ through the strategy as that for the zero-range process given in Section 3.4. below. Note that for the exclusion process, we do not need the replacement lemma, nor the entropy estimate. Then, applying Proposition \ref{be} with $ A=I_d$, we obtain the hydrodynamic limit for $ \Phi $ and Theorem \ref{exclusion} is proved.

\vspace{0.3cm}
\subsection{Zero range process on crystal lattices}

Let $X_N:=(V_N,E_N)$ be the $N$-scaling finite graph. Let $Z_N= \mathbb{N} ^ {V_N}$ be the configuration space and $\eta=\{\eta_x\}_{x\in V_N} \in Z_N$ be the configuration. 
For each $\eta \in Z_N$, $\eta^e$ is obtained from $\eta$ where one particle jumped from $oe$ to $te$, i.e.
\begin{equation*}
\eta^e_x=\left\{
\begin{aligned}
&\eta_{oe}-1\ \ \ \ \ \  x=oe
\\
&\eta_{te}+1\ \ \ \ \ \ x=te
\\
&\eta_x\ \ \ \ \ \ \ \ \ \ \ \  otherwise. 
\end{aligned}
\right.
\end{equation*}

Let $g: \mathbb{N} \to \mathbb{R}_+$ be a function with $g(0)=0$. We assume that $g(k)>0$ for $k >0$ and satisfies the condition
\begin{equation*}
g^*:= \sup_{k \ge 0} | g(k+1)-g(k) |< \infty. 
\end{equation*}
Define 

$$Z(\varphi):=\sum_{k \ge 0} \frac{\varphi^k}{g(k)!},$$
where $g(k)!:=g(1)\cdots g(k), g(0)!:=1$ and let $\varphi^*:=\{\lim\sup_{k\to \infty} \sqrt [k]{g(k)!}\}^{-1}$ be the radius of convergence. Furthermore, we assume that

\begin{equation*}
\lim_{\varphi \uparrow \varphi^*} Z(\varphi)=\infty.
\end{equation*}
Since we need to calculate the exponential moments in our proofs, for the existence of all exponential moments, we suppose $Z(\cdot)$ is finite on $\mathbb{R}_+$, equivalently, $\varphi^*=\infty$. 

For each $ \varphi \ge 0 $, let $\bar{\nu}_\varphi$ be the product measure on $Z_N$ with the marginals given by

\begin{equation}
\bar{\nu}_\varphi \{ \eta; \eta_x=k\} =\frac{1}{Z(\varphi)} \frac{\varphi^k}{g(k)!}.
\end{equation}

 The generator acting on $L^2(Z_N, \bar{\nu}_\varphi)$ as 
\begin{equation}
L_Nf(\eta):=\sum_{e \in E_N} p(e) g(\eta_{oe}) [f(\eta^e) -f(\eta)],\ \ \ \ \ \ \ f \in L^2(Z_N, \bar{\nu}_\varphi)
\end{equation}

\noindent defines a Markov process $\eta(t)$ on $Z_N$ called zero range process with parameters $(p,g)$. 

\begin{proposition}
For each $ \varphi \ge 0$, the product measure $\bar{\nu}_\varphi$ is invariant for the zero range process with parameter $(p,g)$.
\end{proposition}




 
  Let $R(\varphi)=E_{\bar{\nu}_\varphi}[\eta_{x}]$ denote the expection of particles per site. Notice that $R(\varphi)$ is strictly increasing with $R(0)=0$. Denote by $\Psi$ the inverse mapping of $R$. For each $\alpha \ge 0$, define the product measure $\nu_\alpha$ by

\begin{equation}
\nu_\alpha(\cdot)=\bar{\nu}_{\Psi(\alpha)}(\cdot).
\end{equation} 

\noindent Then we obtained a family of invariant measures $\{ \nu_\alpha, \alpha \ge 0\}$ parametrized by density, i.e., for every $\alpha\ge0$,

\begin{equation}
E_{\nu_\alpha}[\eta_x]=\alpha.
\end{equation}
Moreover, it is not hard to see that

\begin{equation}
\Psi(\alpha)=E_{\nu_\alpha} [g(\eta_x)].
\end{equation}
Furthermore, the function $\Psi$ is uniformly Lipschitz on $\mathbb{R}_+$ with Lipschitz constant $g^*$.(See \cite{KL1999})

For arbitrary fixed time $T>0$, let $D([0,T],Z_N) $ be the path space and for a probability measure $\mu^N$ on $Z_N$, let $P_{\mu^N}$ be the distribution on $D([0,T],Z_N) $ of the continuous Markov process with generator $N^2L_N$ and initial measure $\mu^N$ and $E_{\mu^N}$ be the expectation with respect to $P_{\mu^N}$.

\subsection{Relative entropy and Dirichlet form.}

Let $\mathcal{P}(Z_N)$ be the space of probability measures on $Z_N$. For an invariant measure $\nu_{\alpha^*}$ with $\alpha^*>0$, the relative entropy and Dirichlet form of $\mu \in \mathcal{P} (Z_N)$, with respect to $\nu_{\alpha^*}$ is defined as the following:

\begin{equation}
H_N( \mu | \nu_{\alpha^*}):=\int f \log f d\nu_{\alpha^*}
\end{equation}

\begin{equation}
D_N(\mu):=- \langle \sqrt{f}, L_N \sqrt{f} \rangle _{\nu_{\alpha^*}}
\end{equation}
where $f:=\frac{d \mu} {d \nu_{\alpha^*}}$, i.e. $f(\eta)=\frac{ \mu(\eta)} {\nu_{\alpha^*} (\eta)}$. To keep notation simple, we shall denote the entropy by $H_N(f)$ and the Dirichlet form by $D_N(f)$.



\noindent For the zero range process, we have the following result:

\begin{theorem} \label{main}
Let $\Phi$ be a periodic realization with lattice group $\phi(\Gamma)$ and $\Phi_h:=\Phi_\phi^h$ be the harmonic realization associated to the lattice group $\phi(\Gamma)$.
Let $\rho_0 : \mathbb{T}^d_\phi \to \mathbb{R}_+$ be an integrable function. Assume that there exist some positive constants $K_0, \alpha^*$ such that 
\begin{equation}
H_N(\mu^N  | \nu_{\alpha^*}) \le K_0 |V_N|,
\end{equation}


\begin{equation} \label{associate}
\varlimsup_{N \to \infty} \mu^N\left[ \left|\frac{1}{|V_N|} \sum_{x \in V_N} G(\Phi_N(x)) \eta_x -\int_{\mathbb{T}_\phi^d} G(u) \rho_0(u) \frac{du}{vol({\mathbb{T}^d_\phi})}\right|>\delta\right]=0,
\end{equation}
for every $\delta>0$ and every continuous function $G: \mathbb{T}^d_\phi \to \mathbb{R}$.  Then, for every $t>0$,


\begin{equation}
\varlimsup_{N \to \infty} P_{\mu^N } \left[ \left|\frac{1}{|V_N|} \sum_{x \in V_N} G(\Phi_N(x)) \eta_x(t) -\int_{\mathbb{T}_\phi^d} G(u) \rho(t,u) \frac{du}{vol(\mathbb{T}^d_\phi)}\right|>\delta\right]=0,
\end{equation}
for every $\delta>0$ and every continuous function $G: \mathbb{T}^d_\phi \to \mathbb{R}$, where $\rho(t,u)$ is assumed to be the unique weak solution of the following equation

 \begin{equation} \label{eq4}
 \left\{
 \begin{aligned}
  &\frac{\partial }{\partial t} \rho =\nabla \mathbb{D}_{\Phi_h} \nabla  \Psi(\rho),
  \\
 & \rho(0,\cdot)=\rho_0(\cdot).
  \end{aligned}
  \right.
  \end{equation}

\end{theorem}


\begin{remark} \label{graph}
Though the processes on the square lattice are categorized into nearest neighbor interaction models and more general finite range interaction models, the processes on the crystal lattice are assumed to have nearest neighbor interaction without loss of  generality. In fact, for a given set of vertices $V$ and symmetric jump rate  $p(x,y)_{x,y \in V}$, we define $E:=\{(x,y); p(x,y)=p(y,x) >0\}$, then the jumps occur only between nearest neighbor sties by definition and moreover, $p(e)>0$ for all $e \in E$.   In particular, any finite range simple exclusion process (zero range process) on the square lattice can be regarded as a (nearest neighbor) simple exclusion process (zero range process) on a crystal lattice.
\end{remark}

\begin{remark} \label{uniqueness}
We have assumed the uniqueness of the weak solution of equation (\ref{eq4}) with initial condition $ \rho _0 $. In the case of square lattice $ \mathbb {Z}^d$, results on uniqueness can be found in Chapter 5 of \cite{KL1999}. Proceeding as in \cite{KL1999}, we can prove that all limit points of the sequence $\{ \pi_\cdot^{ \Phi, N} \} $ are concentrated on paths satisfying an energy estimate. If $|V_0| =1$, we can obtain the uniqueness through the $L^2$ estimate as in \cite{KL1999}. But more generally($| V_0 | \ge 2 $) it is open. 
\end{remark}


\subsection{Proof of Theorem \ref{main}}
In this subsection, we give the proof of Theorem \ref{main} for the harmonic realization $\Phi$. Then,  with the same discussion mentioned in Proposition \ref{be}, we obtain the hydrodynamic limit of the zero range process for general periodic realizations. And we omit the extension from the harmonic realization to the general periodic realization.  Fix a harmonic realization $\Phi$ and recall that the empirical density is given by
\begin{equation}
\pi_t^N(du)=\pi_t^{\Phi,N}(du):=\frac{1}{|V_N|} \sum_{x \in V_N} \eta_x(t) \delta_{\Phi_N(x)} (du)
\end{equation}

\noindent where $\delta_x$ is the Dirac measure at $x \in \mathbb{T}^d_\phi$, then $\pi_t^N$ is a process taking value in $\mathcal{M}_+$, the space of all finite positive measure on $\mathbb{T}^d _\phi$. Let $D([0,T], \mathcal{M}_+)$ be the path space and $Q^N$ be the measure on $D([0,T], \mathcal{M}_+)$ associated to $\pi^N_t$ starting from $\mu^N$. For every $G \in C^2(\mathbb{T}^d _\phi )$, we define

\begin{equation}
\langle \pi^N_t, G \rangle:=\frac{1}{|V_N|} \sum_{x \in {V_N}} \eta_x(t) G(\Phi_N(x))
\end{equation}
We will consider the martingale $M_t$ and its quadratic process $N_t$ defined by
\begin{equation}
 M_t:=\langle \pi^N_t, G \rangle-\langle \pi^N_0, G \rangle- \int_0^t N^2 L_N \langle \pi^N_s, G \rangle ds
\end{equation}

\begin{equation}
 N_t:=M^2_t - \int_0^t N^2 \left \{ L_N\langle \pi^N_s, G \rangle ^2- 2\langle \pi^N_s, G \rangle L_N \langle \pi^N_s, G \rangle \right\} ds
\end{equation}
Note that $M_t$ and $N_t$ can be rewritten as

\begin{equation}
 M_t:=\langle \pi^N_t, G \rangle-\langle \pi^N_0, G \rangle- \int_0^t \frac{1}{|V_N|} \sum_{x \in V_N} \Delta_NG(\Phi_N(x)) g(\eta_x) ds
\end{equation}

\begin{equation}
 N_t:=M^2_t - \int_0^t  \frac{N^2}{|V_N|^2} \sum_{e \in E_N} p(e) g(\eta_{oe}) (G(\Phi_N(te))-G(\Phi_N(oe)))^2ds
\end{equation}

\noindent where the discrete Laplacian is defined by
\begin{equation}
\Delta_NG(\Phi_N(x)):=2N^2 \sum_{e \in E_x} p(e) [G(\Phi_N(te)) -G(\Phi_N(oe))]
\end{equation}
Since $\Phi$ is harmonic and by Taylor's formula, we have that

\begin{align*} 
\Delta_NG&(\Phi_N(x)):=2N^2 \sum_{e \in E_x} p(e) [G(\Phi_N(te)) -G(\Phi_N(oe))]
\\
&= \sum_{e \in E_x} p(e)\left \{ 2N \sum_{i=1}^d \frac{ \partial G}{ \partial x_i} (\Phi_N (x)) v_i(e)+ \sum_{i,j=1}^d \frac{ \partial^2 G}{ \partial x_i\partial x_j} (\Phi_N(x)) v_i(e)v_j(e) \right\} +o_N
\\
&=\sum_{e \in E_x} \sum_{i,j=1}^d \frac{ \partial^2 G}{ \partial x_i\partial x_j} (\Phi_N(x)) p(e)v_i(e)v_j(e)+o_N
\\
&=F_N(x)+o_N
\end{align*}

\noindent where 
\begin{equation}
F_N(x):=\sum_{e \in E_x} \sum_{i,j=1}^d \frac{ \partial^2 G}{ \partial x_i\partial x_j} (\Phi_N(x)) p(e)v_i(e)v_j(e).
\end{equation}

\vspace{0.3cm}
\begin{lemma} \label{doob}
\begin{equation}
\lim_{N \to \infty} E_{\mu^N} \left [ \sup_{0\le t \le T} |M_t|^2\right]=0
\end{equation}
\end{lemma}
\noindent Proof: First, by entropy inequality, we have that

\begin{equation*}
\begin{aligned}
E_{\mu^N} \left [ \frac{1}{|V_N|} \sum_{x \in V_N} \eta_x\right] &\le \frac{1}{\gamma |V_N|} \log E_{\nu_{\alpha^*}} \left[e^{\gamma \sum_{x \in V_N } \eta _x}\right]+\frac{H_N(\mu^N | \nu_{\alpha^*})}{\gamma |V_N|}\\
&\le \frac{1}{\gamma} \log E_{\nu_{\alpha^*}} \left[e^{\gamma  \eta _x}\right]+\frac{K_0}{\gamma}
\end{aligned}
\end{equation*}
for all $\gamma>0$. Furthermore, by Dood's inequality, 

\begin{align*}
E_{\mu^N} \left [ \sup_{0\le t \le T} |M_t|^2\right] &\le 4E_{\mu^N} \left [ |M_T|^2 \right]
\\
&=E_{\mu^N} \left [  \int_0^T \frac{N^2}{|V_N|^2} \sum_{e \in E_N} p(e) g(\eta_{oe}) (G(\Phi_N(te))-G(\Phi_N(oe)))^2ds\right]
\\
&\le \frac{4 C( p, G ) g^* T}{|V_N|} E_{\mu^N} \left [ \frac{1}{|V_N|} \sum_{x \in V_N} \eta_x\right] \to 0
\end{align*}
as $N \to \infty$, where $C( p, G )$ is a constant depending only on $p( \cdot )$ and $G( \cdot )$. 

\vspace{0.3cm}
\begin{lemma} \label{cpt}
The sequence of $\{Q^N\}$ is relatively compact and all limit points $Q^*$ are concentrated on trajectories of absolutely continuous measures with respect to the Lebesgue measure:

\begin{equation}
Q^* [\pi : \pi_t(du)=\pi(t,u) du]=1
\end{equation}
\end{lemma}

\noindent Proof: For each $G \in C^2(\mathbb{T}^d_\phi)$, define the probability measure $Q^NG^{-1}$ on $D([0,T],\mathbb{R})$ by setting
\begin{equation*}
Q^NG^{-1}(\cdot)=Q^N[\{\pi\ |\ \langle \pi, G \rangle \in \cdot\}]
\end{equation*}
To show $\{Q^N\}$ is relatively compact, it suffices to show that, for each $G \in C^2(\mathbb{T}^d_\phi)$, the following (i) and (ii) holds.(See \cite{KL1999})
\\
(i) For any $t>0$ and any $\epsilon >0$, there exists a constant $A>0$ such that 
\begin{equation*}
\sup_{N}Q^NG^{-1}[|\pi_t|>A]\le\epsilon
\end{equation*}
Since 

\begin{align*}
\sup_{N}Q^NG^{-1}[|\pi_t|>A]&=\sup_NQ^N[\{\pi\ |\ |\langle \pi_t, G \rangle| >A\}]
\\
&=\sup_NP_{\mu^N}[|\langle \pi_t^N, G \rangle| >A\}]
\\
&=\sup_N P_{\mu^N} \left[ \left| \frac{1}{|V_N|} \sum_{x \in V_N} G(\Phi_N(x)) \eta_x(t) \right|>A\right]
\\
&\le \sup_N P_{\mu^N} \left [ \frac{1}{|V_N|} \sum_{x \in V_N} \eta_x(t) > \frac{A}{||G||_\infty}\right]
\\
&=\sup_N \mu^N \left [ \frac{1}{|V_N|} \sum_{x \in V_N} \eta_x > \frac{A}{||G||_\infty}\right]
\end{align*}
In the last equation, we used the conservation of the number of total particles. By (\ref{associate}), for any $ \epsilon >0$, there exists $N_0>0$ such that
\begin{equation*}
\sup_{N\ge N_0} \mu^N \left [ \frac{1}{|V_N|} \sum_{x \in V_N} \eta_x > ||\rho_0||_{L^1}+1\right] \le \epsilon.
\end{equation*}
Take $A$ large enough such that 
$$\frac{A}{||G||_\infty}\ge ||\rho_0 ||_{L^1}+1,$$
$$\sup_{N< N_0} \mu^N \left [ \frac{1}{|V_N|} \sum_{x \in V_N} \eta_x >\frac{A}{||G||_\infty} \right] \le \epsilon.$$
For such $A$, we have that
$$\sup_N \mu^N \left [ \frac{1}{|V_N|} \sum_{x \in V_N} \eta_x > \frac{A}{||G||_\infty}\right]\le \epsilon.$$
Thus, (i) holds.
\vspace{0.2cm}

\noindent(ii) For any $\epsilon>0$, it holds that
\begin{equation*}
\lim_{\gamma \to 0} \varlimsup_{N \to \infty} Q^NG^{-1} \left[\pi\ |\ \sup_{|s-t|<\gamma} |\pi_t-\pi_s|>\epsilon \right]=0.
\end{equation*}
Note that 
$$ Q^NG^{-1} \left [\pi\ |\ \sup_{|s-t|<\gamma} |\pi_t-\pi_s|>\epsilon \right ]=P_{\mu^N}  \left [ \sup_{|s-t|<\gamma} | \langle \pi_t^N,G \rangle -\langle \pi_s^N,G \rangle |>\epsilon \right],$$
and 
$$|\langle \pi_t^N,G \rangle -\langle \pi_s^N,G \rangle | \le | M_t-M_s | +\int_s^t \left |\frac{1}{|V_N|} \sum_{x \in V_N} \Delta_N G(\Phi_N(x)) g(\eta_x)\right | dr.$$
Furthermore, we have 

\begin{equation*}
\begin{aligned}
\lim_{\gamma \to 0}\varlimsup_{N \to \infty} P_{\mu^N} &\left [ \sup_{|s-t|<\gamma} \int_s^t \left |\frac{1}{|V_N|} \sum_{x \in V_N} \Delta_N G(\Phi_N(x)) g(\eta_x)\right | dr>\epsilon \right]
\\
& \le \lim_{\gamma \to 0}\varlimsup_{N \to \infty} \mu^N \left [ C( p, G ) g^* \gamma \frac{1}{|V_N|} \sum_{x \in V_N} \eta_x >\epsilon\right ]=0.
\end{aligned}.
\end{equation*}
By Chebyshev's inequality, Doob's inequality and Lemma \ref{doob}, we have
\begin{equation*}
\lim_{\gamma \to 0}\varlimsup_{N \to \infty} P_{\mu^N} \left [ \sup_{|s-t|<\gamma}  |M_t-M_s|>\epsilon \right] \le \lim_{\gamma \to 0}\varlimsup_{N \to \infty} \frac{4}{\epsilon^2}E_{\mu^N} [M^2_T] =0.
\end{equation*}
Then, by triangular inequality, (ii) holds. Thus, $\{Q^N\}$ is relatively compact. For the left part and more details of the proof, see Section 5.1 of \cite{KL1999}.

\vspace{0.3cm}
We now prove Theorem \ref{main} by assuming the following replacement lemma and give its proof in Section 4. Before giving a statement of this lemma, we introduce two local averages. 

Take a $\mathbb{Z}$-basis $\{ \sigma_1, \dots, \sigma_d\}$ of $\Gamma$ and identify $\Gamma$ with $\mathbb{Z}^d$. We define the standard generator system of $\Gamma$ by setting
      \begin{equation*}
      S=\{  \sigma_1, \dots, \sigma_d,  -\sigma_1, \dots, -\sigma_d\}.
      \end{equation*}
      
     \noindent We introduce the length function associated to $S$, $| \cdot |: \Gamma \to \mathbb{N}$ by
     \begin{equation*}
     | \sigma |:=\min \left \{\ l\ |\ \sigma =\sum_{k=1}^l \epsilon _{i_k} \sigma _{i_k}, \epsilon_{i_k} \in \{-1,1\}, i_k \in \{ 1,\cdots,d\} \right \},
     \end{equation*}
for $\sigma \in \Gamma$. Then the map $(\sigma ,\sigma ^\prime) \in \Gamma \times \Gamma \to |\sigma -\sigma^\prime | \in \mathbb{N}$ induces a metric in $\Gamma$, which is called the word metric 
associated to $S$. By the natural homomorphism from $\Gamma$ to $\Gamma_N$, the length function and metric are defined in the same way. To abuse the notation, we use the same symbol. Since $\Gamma$ acts on $X$ freely, for each $x \in V$, there exists a unique $\sigma_x \in \Gamma$ such that $x \in \sigma_x D_{x_0}$. Define the map $[\cdot]: V \to \Gamma$ by setting $[x]=\sigma_x$ for $x \in V$.

For $R>0$, define the $R-$ball by setting

\begin{equation}
B(D_{x_0},R):=\bigcup_{\sigma \in \Gamma, |\sigma|\le R} \sigma D_{x_0}.
\end{equation}
Define a local average of $\eta$ for $x \in V$ by

\begin{equation}
\bar{\eta}_{x,R}:=\frac{1}{| [x] B(D_{x_0}, R)|} \sum_{z \in [x]B(D_{x_0},R)}\eta_z
\end{equation}
and a local average of $g$ by
\begin{equation}
\widetilde {g}_{x,\epsilon N}(\eta):=\frac{1} {|\{ \underline{\sigma} | |\underline {\sigma}| \le \epsilon N\}|} \sum_{ |\underline {\sigma} | \le \epsilon N} g_{\underline{\sigma} x} (\eta)
\end{equation}
where $g_x(\eta):=g(\eta_x)$ for $x \in V$ and $ | U |$ stands for the number of the elements in $U$. Note that $ \bar { \eta }_{ x , R }= \bar { \eta}_{ x_0, R } $ for every $ x \in D_{ x_0 } $.

\begin{lemma} \label{replace} {(Replacement lemma)} 
For every $x \in D_{x_0}$, 
\begin{equation} \label {replaceeq}
\varlimsup_{\epsilon \to 0} \varlimsup_{N  \to \infty} P_{\mu^N}\left[\int _0^T \frac{1}{|\Gamma_N|} \sum_{\sigma \in \Gamma_N} V_{\sigma x, \epsilon N} (\eta(t)) dt >\delta\right] =0.
\end{equation}
where 

\begin{equation}
V_{\sigma x, \epsilon N} (\eta):=| \widetilde {g} _{\sigma x, \epsilon N} (\eta) -\Psi (\bar{\eta} _{\sigma x_0, \epsilon N})|. 
\end{equation}

\end{lemma}

\noindent Proof of Theorem \ref{main}: 
For $ G: [0,T] \times \mathbb{T}_{\phi}^d \to \mathbb{R}$ of class $C^{1,2}$, consider the martingale $M_t=M^{G,N}_t$ given by 

\begin{equation*}
 M_t:=\langle \pi^N_t, G_t \rangle-\langle \pi^N_0, G_0 \rangle- \int_0^t ( \partial_s +N^2L_N) \langle \pi^N_s, G_s \rangle ds
\end{equation*}
As the above Lemma \ref{doob}, we have that 




\begin{equation*}
\lim_{N\to \infty} P_{\mu^N} \left[ \sup_{0 \le t \le T} \left | \langle \pi^N_t, G_t \rangle-\langle \pi^N_0, G_0 \rangle-\int_0^t   \langle \pi^N_s, \partial_sG_s \rangle ds- \int_0^t \frac{1}{|V_N|} \sum_{x \in V_N} \Delta_NG(\Phi_N(x)) g(\eta_x) ds \right |>\delta
\right] =0
\end{equation*}
Note that $g(\eta_x)$ is not a function of the empirical process $\pi^N_{\cdot}$, to close the equation, we replace $g(\eta_x)$ in the following steps.

\vspace{0.2cm}

\noindent(i) Firstly, we can replace $g_x(\eta)$ by $ \widetilde{g}_{x,\epsilon N}(\eta)$.

\vspace{0.2cm}
\noindent Since 
\begin{align*}
P_{\mu^N}  & \left[ \left |\int_0^T \frac{1}{|V_N|} \sum_{x \in V_N} \Delta_N G(\Phi_N(x))  \left \{g_x(\eta(s)) -\widetilde{g} _{x,\epsilon N} (\eta(s)) \right\} ds \right | >\delta \right]
\\
&\le \frac{1}{\delta} E_{\mu^N} \left [\left |\int _0^T \frac{1}{|V_N|} \sum_{x \in D_{x_0}} \sum_{\sigma \in \Gamma_N} g_{\sigma x} (\eta(s))  \left \{ \Delta_N G(\Phi_N(\sigma x)) -\frac{1} {|\{ \underline{\sigma} | |\underline {\sigma}| \le \epsilon N\}|} \sum_{ |\underline {\sigma} | \le \epsilon N} \Delta_N G(\Phi_N( \underline{\sigma} \sigma x))  \right\} \right | \right]
\\
&\le \frac{1}{\delta} E_{\mu_N} \left[\int _0^T g^* \sup_{x \in D_{x_0}\atop \sigma \in \Gamma_N, |\underline{\sigma}| \le \epsilon N} |\Delta_N G(\Phi_N(\sigma x))-\Delta_N G(\Phi_N( \underline{\sigma} \sigma x))| \frac{1}{|V_N|} \sum_{x \in V_N} \eta_{x}(s) ds \right]
\\
&=\frac{g^*T}{\delta} E_{\mu_N} \left [\sup_{x \in D_{x_0} \atop \sigma \in \Gamma_N, |\underline{\sigma}| \le \epsilon N}  \left | \Delta_N G(\Phi_N(\sigma x))-\Delta_N G(\Phi_N( \underline{\sigma} \sigma x)) \right | \frac{1}{|V_N|} \sum_{x \in V_N} \eta_{x}(0) \right ]
\end{align*} 

\noindent which tends to 0 as $N \to \infty, \epsilon \to 0$ by Lebesgue dominated convergence theorem. Note that, we use the conservation of the total number of particles in last equation. 

\vspace{0.2cm}
\noindent(ii) By the replacement lemma, we replace $\widetilde{g} _{x,\epsilon N} (\eta) $ by $\Psi (\bar{\eta}_{ [x] x_0, \epsilon N})$, 

$$\varlimsup_{\epsilon \to 0} \varlimsup_{N \to \infty} P_{\mu^N} \left [ \left |\int_0^T \frac{1}{|V_N|} \sum_{x \in V_N} \Delta_N G(\Phi_N(x)) \{\widetilde{g} _{x,\epsilon N} (\eta(s)) -\Psi (\bar{\eta}_{ [x ] x_0, \epsilon N})\} ds \right | >\delta \right ]=0.$$ 

\noindent(iii) By Taylor's formula, we replace the discrete Laplacian $ \Delta_N G(\Phi_N(x))$ by $F_N(x)$, and obtain that

$$\varlimsup_{\epsilon \to 0} \varlimsup_{N \to \infty} P_{\mu^N} \left [\int _0^T \left | \frac{1}{|V_N|} \sum_{x \in V_N}  \left \{ \Delta_N G(\Phi_N(x))-F_N(x) \right \} \Psi(\bar{\eta}_{[x ] x_0,\epsilon N})\right| ds >\delta \right ]=0.$$



\noindent Note that

\begin{equation*}
F_N(\sigma x)=\sum_{e \in E_x} \sum_{i,j=1}^d \frac{ \partial^2 G}{ \partial x_i\partial x_j} (\Phi_N(\sigma x_0 ) ) p(e)v_i(e)v_j(e)+o_N,\ \ \ \forall x \in D_{x_0},\ \forall \sigma \in \Gamma_N.
\end{equation*}

\noindent Together with (i), (ii) and (iii), we have
\begin{equation*}\small
\begin{aligned}
\varlimsup_{\epsilon \to 0} \varlimsup_{N \to \infty} Q^N \left [ \left | \langle \pi^N_T, G \rangle -\langle \pi^N_0, G \rangle -\int_0^T   \langle \pi^N_s, \partial_sG_s \rangle ds -\int _0^T \frac{1}{|\Gamma_N|} \sum_{\sigma \in \Gamma _N} \nabla  \mathbb{D}_\Phi  \nabla G(\Phi_N(\sigma x_0)) \Psi(\bar{\eta}_{\sigma x_0, \epsilon N}) ) ds \right | > \delta \right]=0
\end{aligned}
\end{equation*}

\noindent By Lemma \ref{chi}, we can replace $\Psi(\bar{\eta}_{\sigma x_0, \epsilon N}) )$ by $\Psi(\langle \pi_t^N, \chi _{\Phi_N(\sigma x_0), \epsilon} \rangle)$, i.e.,

\begin{equation*}\small
\begin{aligned}
\varlimsup_{\epsilon \to 0} \varlimsup_{N \to \infty} Q^N \left [ \left | \langle \pi^N_T, G \rangle -\langle \pi^N_0, G \rangle  - \int_0^T   \langle \pi^N_s, \partial_sG_s \rangle ds -\int _0^T \frac{1}{|\Gamma_N|} \sum_{\sigma \in \Gamma _N} \nabla  \mathbb{D}_\Phi  \nabla G(\Phi_N(\sigma x_0)) \Psi( \langle \pi_s^N, \chi _{\Phi_N(\sigma x_0), \epsilon} \rangle) ds \right | > \delta \right ]=0
\end{aligned}
\end{equation*}

\noindent By Lemma \ref{cpt}, for any limit point $Q^*$ of $\{Q^N\}$, we have that 

$$\varlimsup_{\epsilon \to 0} Q^* \left [ \left | \langle \rho_T, G \rangle -\langle \rho_0, G \rangle -\int _0^T \int _{\mathbb{T}^d_\phi} \left\{ \rho(s,u)\partial_sG(u)+\nabla \mathbb{D}_\Phi \nabla G(u)  \cdot \Psi (\langle \rho_s, \chi_{\cdot,\epsilon} \rangle ) \frac{du}{vol(\mathbb{T}^d_\phi)} ds \right\} \right | >\delta \right]=0$$ 

\noindent By the dominated convergence theorem, as $\epsilon \to 0$, we obtain that,

$$ Q^* \left [ \left | \langle \rho_T, G \rangle -\langle \rho_0, G \rangle -\int _0^T \left\{  \rho(s,u)\partial_sG(u)+  \langle \Psi (\rho_t ),  \nabla \mathbb{D}_\Phi \nabla G \rangle dt \right\}  \right | >\delta \right ]=0$$

\noindent Up to now, we have proved that $Q^*$ concentrates on path $\{ \rho_t, 0\le t \le T\}$, which is the weak solution of (\ref{eq1}). In particular, we assume such a weak solution is unique, then we have that, for every $\delta>0$, 
\begin{equation*}
\lim_{N \to \infty} Q^N \left [ d_{S_k} \left (\pi^N,  \frac{\rho(\cdot,u)du}{vol(\mathbb{T}^d_\phi)} \right) >\delta \right]=0
\end{equation*}

\noindent where $d_{S_k}$ is the Skorohod distance on $D([0,T], \mathcal{M}_+)$. Since the limit measure is concentrated on weakly continuous trajectories, for fixed time $t \in (0,T)$, $\pi^N_t$ converges in distribution to the deterministic measure $\rho(t,u) du/ vol(\mathbb{T}^d_\phi)$. Since convergence in distribution to a deterministic variable implies convergence in probability. In particular, since $T>0$ is arbitrary, we obtain that for every $\delta>0$ ,every $t>0$ and every continuous function $G \in C(\mathbb{T}^d_\phi)$, 

$$\varlimsup_{N \to \infty} P_{\mu^N } \left [  \left |\frac{1}{|V_N|} \sum_{x \in V_N} G(\Phi_N(x)) \eta_x(t) -\int_{\mathbb{T}^d_\phi} G(u) \Psi (\rho(t,u) ) \frac{du}{vol(\mathbb{T}^d_\phi)} \right |>\delta \right ]=0$$
which concludes the Theorem 3.2.

\section{Replacement Lemma}
In this section, we prove the replacement lemma through one block estimate and two blocks estimate. 

\subsection{Proof of replacement lemma}
Let $\mu^N_t$ be the distribution of $\eta(t)$ on $Z_N$ and set 

\begin{equation}
f_t^N:=\frac{ d \mu^N_t} { d \nu_{\alpha^*}},\ \ \ \ \ \ \ \ \ \ \ \bar{f}_T^N:= \frac{1}{T} \int_0^T f^N_t dt .
\end{equation}
Then $f_t^N$ satisfies the Kolmogorov equation

\begin{align*}
\left \{
\begin{array}{lll}
\partial_t f_t^N&=N^2 L_N f_t^N,\\
f^N_0&=\frac{d \mu^N} {d \nu_{\alpha^*}}.
\end{array}
\right.
\end{align*}
\noindent As in Section 5.2 of \cite{KL1999}, we have that


\begin{equation} \label{Dirichlet}
H_N(\bar{f}^N_T) \le H_N(f^N_0),\ \ \ \ \ \ \ \ \ \ \ D_N(\bar{f}^N_T) \le \frac {H_N(f^N_0)} {2T N^2}.
\end{equation}
Note that
\begin{align*}
P_{\mu^N}& \left[ \int _0^T \frac{1}{|\Gamma_N|} \sum_{\sigma \in \Gamma_N} V_{\sigma x, \epsilon N} (\eta(t)) dt >\delta \right]  
\\
&\le \frac {1} {\delta} E_{\mu^N} \left [\int _0^T \frac {1} { | \Gamma_N |} \sum_{\sigma  \in \Gamma_N} V_{\sigma x, \epsilon N} (\eta(t)) dt \right]
\\
&= \frac {1} {\delta} \int _0^T \int \frac {1} {| \Gamma_N |} \sum_{\sigma  \in \Gamma_N} V_{\sigma x, \epsilon N} (\eta) f^N_t(\eta) \nu_{\alpha^*}(d \eta)dt
\\
&=\frac {T} {\delta} \int \frac {1} { | \Gamma_N |} \sum_{\sigma  \in \Gamma_N} V_{\sigma x, \epsilon N} (\eta) \bar{f}^N_T(\eta) \nu_{\alpha^*}(d \eta).
\end{align*}
Thus, to show the replacement lemma, it suffices to show the following lemma.

\begin{lemma} \label{re2}
 For every $C>0$ and every $x \in D_{x_0}$,
\begin{equation}
\varlimsup_{\epsilon \to 0} \varlimsup_{ N \to \infty} \sup_{H_N(f) \le C|V_N|, \atop
D_N(f) \le C\frac{|V_N|}{N^2}}\int \frac {1} {| \Gamma_N |} \sum_{\sigma  \in \Gamma_N} V_{\sigma x, \epsilon N} (\eta) f(\eta) \nu_{\alpha^*}(d \eta)=0.
\end{equation}

\end{lemma}
\noindent Lemma \ref{re2} will be shown by the following one block estimate and two blocks estimate.

\begin{lemma}(One block estimate) \label{one}
For every $C>0$ and every $x \in D_{x_0}$,

\begin{equation} 
\begin{aligned}
\varlimsup_{l \to \infty} \varlimsup_{ N \to \infty} \sup_{H_N(f) \le C|V_N| \atop D_N(f) \le C\frac{|V_N|}{N^2}}\int \frac {1} { | \Gamma_N |} \sum_{\sigma  \in \Gamma_N} V_{\sigma x, l} (\eta) f(\eta) \nu_{\alpha^*}(d \eta)=0.
\end{aligned}
\end{equation}
\end{lemma}

\begin{lemma}(Two blocks estimate) \label{two}
For every $C>0$ and every $x \in D_{x_0}$, 
\begin{equation}
\begin{aligned}
\varlimsup_{l \to \infty} \varlimsup_{\epsilon \to 0} \varlimsup_{N \to \infty} \sup_{H_N(f) \le C|V_N| \atop D_N(f) \le C \frac {|V_N|} {N^2}} \sup_{|\underline{\sigma}| \le \epsilon N} \int \frac {1} { | \Gamma_N |} \sum_{\sigma \in \Gamma_N} |\bar{\eta} _{\underline{\sigma} \sigma x,l} -\bar{\eta}_{\sigma x, \epsilon N}| f(\eta) \nu_{\alpha^*} (d \eta)
\end{aligned}
\end{equation}
\end{lemma}

We first prove Lemma \ref{replace} by using the one block estimate and two blocks estimate.
Recall that

$$V_{\sigma x, \epsilon N}(\eta) =\left | \frac{1}{|\{ \underline {\sigma} | |\underline{\sigma}| \le \epsilon N\}|} \sum_{|\underline {\sigma}| \le \epsilon N} g_{ \underline {\sigma} \sigma x} -\Psi(\bar{\eta} _{\sigma x_0, \epsilon N}) \right | $$
Add and subtract the following expression 

$$\frac{1}{|\{ \underline {\sigma} | |\underline{\sigma}| \le \epsilon N\}|} \sum_{|\underline {\sigma}| \le \epsilon N}  \frac{1}{|\{  \sigma^\prime | |\sigma^\prime| \le l\}|} \sum_{|\sigma^\prime| \le l} \{ g_{\sigma^\prime  \underline {\sigma} \sigma x} (\eta) -\Psi (\bar{\eta} _{ \underline {\sigma} \sigma x_0 , l}) \}$$
By the triangular inequality, we can estimate the integral term appeared in (\ref{replaceeq}) by three parts separately.

$$\int \frac {1} { | \Gamma_N |} \sum_{\sigma  \in \Gamma_N} V_{\sigma x, \epsilon N} (\eta) f(\eta) \nu_{\alpha^*}(d \eta) \le I_1+I_2+I_3$$ 
where 

\begin{align*} 
I_1&:= \int \frac {1} {|\Gamma_N|} \sum_{\sigma \in \Gamma_N}  \left | \widetilde {g}_{\sigma x, \epsilon N} (\eta) -\frac{1}{|\{ \underline {\sigma} | |\underline{\sigma}| \le \epsilon N\}|} \sum_{|\underline {\sigma}| \le \epsilon N}  \frac{1}{|\{  \sigma^\prime | |\sigma^\prime| \le l\}|} \sum_{|\sigma^\prime| \le l}  g_{\sigma^\prime \underline {\sigma} \sigma x} (\eta)
 \right | f(\eta) \nu_{\alpha^*} (d \eta)\\
&=\int \frac {1} {|\Gamma_N|} \sum_{\sigma \in \Gamma_N}  \left | \frac{1}{|\{ \underline {\sigma} | |\underline{\sigma}| \le \epsilon N\}|} \sum_{|\underline {\sigma}| \le \epsilon N} \left[g_{\underline{\sigma} \sigma x} -\frac{1}{|\{  \sigma^\prime | |\sigma^\prime| \le l\}|} \sum_{|\sigma^\prime| \le l}  g_{\sigma^\prime \underline {\sigma} \sigma x} (\eta) \right]
 \right | f(\eta) \nu_{\alpha^*} (d \eta)\\
& \le \int \frac {1} { |\Gamma_N|} \sum_{\sigma \in \Gamma_N} \frac{1}{ (2 \epsilon N +1)^{d}} \sum_{\epsilon N -l \le \underline {\sigma } \le \epsilon N +l} (2l+1)^d g_{\underline {\sigma} \sigma x} (\eta) f(\eta) \nu_{\alpha^*} (d\eta)\\
&\le \frac{C_1(l,g^*,d)}{ \epsilon N} \int \frac{1}{ | \Gamma_N |} \sum_{\sigma \in \Gamma_N} \eta_{\sigma x} f(\eta) \nu_{\alpha^*} (d\eta)
\end{align*}
where $C_1$ is a constant depending on $l,d$ and $g^*$. By the entropy inequality, $I_1$ is bounded above by
\begin{align*}
I_1 &\le\frac{C_1(l,g^*,d)}{\epsilon N |\Gamma_N| \gamma} \left \{ H_N(f) +\log E_{\nu_{\alpha^*}} \left [e^{\sum_{\sigma \in \Gamma_N} \gamma \eta _{\sigma x}}\right ] \right\}
\\
&\le \frac{C_1(l,g^*,d)}{\epsilon N |\Gamma_N| \gamma} \left \{ C|V_N| + |\Gamma_N| \log E_{\nu_{\alpha^*}} \left [ e^{\gamma \eta _{\sigma x}} \right ] \right \}
\\
&=\frac{C_1(l,g^*,d)}{\epsilon N \gamma} \left \{ C|D_{x_0}| +\log E_{\nu_{\alpha^*}} \left [ e^{ \gamma \eta _{\sigma x}} \right ] \right \}
\end{align*}
which tends to 0 as $N \to \infty$.

\begin{align*}
I_2&:= \int \frac {1} {| \Gamma_N |} \sum _{\sigma \in \Gamma_N} \left | \Psi (\bar {\eta} _{\sigma x_0, \epsilon N}) -\frac{1}{|\{ \underline {\sigma} | |\underline{\sigma}| \le \epsilon N\}|} \sum_{|\underline {\sigma}| \le \epsilon N}  \Psi (\bar{\eta} _{ \underline {\sigma} \sigma x_0 , l}) \right | f(\eta) \nu_{\alpha^*} (d \eta)
\\
&\le \sup_{ | \underline {\sigma}| \le \epsilon N} \int \frac{1} { |\Gamma_N|} \sum_{\sigma \in \Gamma_N} \left |\Psi (\bar {\eta} _{\sigma x_0, \epsilon N})-\Psi (\bar {\eta} _{ \underline {\sigma}\sigma x_0, l}) \right | f(\eta) \nu_{\alpha^*} (d\eta)
\\
& \le g^* \sup_{ | \underline {\sigma}| \le \epsilon N} \int \frac{1} { |\Gamma_N|} \sum_{\sigma \in \Gamma_N} |\bar {\eta} _{\sigma x_0, \epsilon N}-\bar {\eta} _{ \underline {\sigma}\sigma x_0, l}| f(\eta) \nu_{\alpha^*} (d\eta)
\end{align*}
which tends to 0 as $N \to \infty, \epsilon \to 0, l \to \infty$ by two blocks estimate.
\begin{align*}
I_3&:= \int \frac{1}{|\Gamma_N|} \sum_{ \sigma \in \Gamma_N} \left | \frac{1}{|\{ \underline {\sigma} | |\underline {\sigma}| \le \epsilon N\}|} \sum_{|\underline {\sigma}| \le \epsilon N}  \frac{1}{|\{  \sigma^\prime | |\sigma^\prime| \le l\}|} \sum_{|\sigma^\prime| \le l} \left\{ g_{\sigma^\prime  \underline {\sigma} \sigma x} (\eta) -\Psi (\bar{\eta} _{ \underline {\sigma} \sigma x_0 , l}) \right\} \right |f(\eta) \nu_{\alpha^*} (d\eta)
\\
&\le \int \frac{1}{ | \Gamma_N |} \sum_{\sigma \in \Gamma_N} |\widetilde {g}_{\sigma x,l} (\eta) -\Psi( \bar{\eta}_{\sigma x_0, l})|f(\eta) \nu_{\alpha^*} (d\eta)
\end{align*}
which tends to 0 as $N \to \infty, l \to \infty$ by one block estimate.
\subsection{Proof of one block estimate} 
Before the proof of Lemma \ref{one}, we prove the following lemma first, which allows us to cut off the large density.

\begin{lemma}
For every $C>0$,
\begin{equation}
\varlimsup_{A \to \infty} \varlimsup_{l \to \infty} \varlimsup_{N \to \infty} \sup_{H_N(f) \le C|V_N|} \int \frac{1} {| \Gamma_N |} \sum_{\sigma \in \Gamma_N} \bar {\eta}_{\sigma x_0,l} 1_{ \{\bar{\eta}_{\sigma x_0,l } \ge A \}} f(\eta) \nu_{\alpha^*} (d\eta)=0
\end{equation}
\end{lemma}
\noindent Proof: By the entropy inequality, the integral part is bounded above by

\begin{equation} \label{cut1}
 \frac{1}{\gamma |\Gamma_N|} \{ H_N(f) +\log E_{\nu_{\alpha*}} [e^{\gamma \sum_{\sigma \in \Gamma_N} \bar {\eta} _{\sigma x_0,l} 1_{ \{ \bar {\eta} _{\sigma x_0,l} \ge A\} }}]\}.
\end{equation}

\noindent Recall that 

$$\bar {\eta} _{\sigma x_0,l}= \frac {1}{ | D_{x_0} | | \{ \underline {\sigma} | |\underline {\sigma}| \le l \} |} \sum_{x \in D_{x_0}} \sum_{| \underline {\sigma} | \le l} \eta _{\underline {\sigma} \sigma x}$$
\\then we have $\bar {\eta} _{\sigma x_0,l}$ and $\bar {\eta} _{\sigma^\prime x_0,l}$ are independent for $ |\sigma -\sigma^\prime| >2l $ under the product measure $ \nu_{ \alpha^*} $. For each $\sigma \in \{ \sigma | |\sigma| \le l\}$, denote $\Omega_\sigma:= \{\underline {\sigma} | \underline {\sigma} -\sigma \in (2l+1) \Gamma \} $. Notice that $ \{ \bar {\eta}_{ \underline {\sigma}  x_0,l} 1_{ \{\bar{\eta}_{\underline {\sigma} x_0,l\} } \ge A} \}_{\underline{\sigma} \in \Omega_{\sigma}} $ are independent. By H\"older's inequality, Chebyshev exponential inequality and $\log (1+x) \le x $ for all $x \ge 0$, we have that

\begin{equation*}
\begin{aligned}
 \frac{1}{\gamma |\Gamma_N|}&\log E_{\nu_{\alpha*}} \left[e^{\gamma \sum_{\sigma \in \Gamma_N} \bar {\eta} _{\sigma x_0,l} 1_{ \{ \bar {\eta} _{\sigma x_0,l} \ge A\} }} \right]
 \\
 &\le  \frac{1}{\gamma |\Gamma_N|}\log \prod_{ | \sigma | \le l} \left \{ E_{\nu_{\alpha*}} \left [e^{\gamma  (2l+1)^d \sum_{ \underline{\sigma} \in \Omega_\sigma} \bar {\eta} _{\underline {\sigma} x_0,l} 1_{ \{ \bar {\eta} _{\underline {\sigma} x_0,l} \ge A\} }} \right ]  \right\} ^{\frac {1} {(2l+1)^d}} 
 \\
 &=\frac{1}{\gamma | \Gamma_N |} \frac{1} {(2l+1)^d} \sum_{ | \sigma | \le l} \sum_{\underline {\sigma} \in \Omega _\sigma } \log E_{\nu_{\alpha*}} \left [e^{\gamma  (2l+1)^d \bar {\eta} _{\underline {\sigma} x_0,l} 1_{ \{ \bar {\eta} _{\underline {\sigma} x_0,l} \ge A\} }} \right] 
 \\
 &\le \frac{1}{ \gamma (2l+1)^d}  \log \left ( 1 + E_{\nu_{\alpha*}} \left [1_{ \{ \bar {\eta} _{ x_0,l} \ge A\} } e^{\gamma  (2l+1)^d \bar {\eta} _{ x_0,l} } \right ] \right )  
 \\
 &\le \frac{1}{ \gamma (2l+1)^d} \left \{   \nu_{\alpha^*}  \left [ \eta | \bar{\eta} _{x_0,l} \ge A \right ] ^{\frac{1}{2}}  \left (E_{\nu_{\alpha*}}  \left [e^{ \frac{2 \gamma} { |D_{x_0} | } \sum_{x \in D_{x_0}} \sum_{ |\sigma| \le l} \eta _{\sigma x}}    \right ] \right )^{\frac{1}{2}} \right \}
\\
&\le \frac{1}{ \gamma (2l+1)^d}  \left \{  e^{-A(2l+1)^d} E_{\nu_{\alpha^*}}  \left [e^{ (2l+1)^d \bar {\eta} _{x_0,l}} \right ] E_{\nu_{\alpha*}}  \left [ e^{ \frac{2 \gamma} { |D_{x_0} | } \sum_{x \in D_{x_0}} \sum_{ |\sigma| \le l} \eta _{\sigma x}}   \right  ]  \right \} ^{\frac{1}{2}}  
\\
&= \frac {1}{ \gamma (2l+1)^d} \left  \{ e^{-\frac{1}{2} (2l+1)^d \left  [A- | D_{x_0} | \log M_{\alpha^*} \left ( \frac{1}{ | D_{x_0} | } \right ) -| D_{x_0} | \log M_{\alpha^*}  \left (\frac{2\gamma} { | D_{x_0} |} \right ) \right ] } \right\}
\end{aligned}
\end{equation*}

\noindent where $M_{\alpha^*} (\theta) :=E_{\alpha^*} [e^{\theta \eta _{x_0}}]$. Then (\ref{cut1}) is bounded by

\begin{equation} \label{cut2}
\frac{C | V_N |} {\gamma | \Gamma_N | } +\frac {1}{ \gamma (2l+1)^d}  \left \{ e^{-\frac{1}{2} (2l+1)^d \left [A- | D_{x_0} | \log M_{\alpha^*}  \left ( \frac{1}{ | D_{x_0} | } \right ) -| D_{x_0} | \log M_{\alpha^*}  
\left (\frac{2\gamma} { | D_{x_0} |} \right ) \right] } \right\}.
\end{equation}
We first fix $\gamma $ large enough so that $\frac{C | V_N |} {\gamma | \Gamma_N | }=\frac{C} {\gamma | D_{x_0} | }$ is small, then pick $A$ large enough such that $A- | D_{x_0} | \log M_{\alpha^*} ( \frac{1}{ | D_{x_0} | }) -| D_{x_0} | \log M_{\alpha^*} (\frac{2\gamma} { | D_{x_0} |})>0$. For such $\gamma$ and $A$, (\ref{cut2}) tends to 0 as $l \to \infty$. The proof is completed.

\vspace{0.5cm}
Now we are ready for the proof of Lemma \ref{one}.

\vspace{0.5cm}
\noindent Proof of Lemma \ref{one}: Notice that

\begin{align*}
V_{ \sigma x, l} &= \left | \frac {1} { | \{ \underline { \sigma } | |\underline {\sigma }| \le l \} |}   \sum_{ | \underline {\sigma } | \le l} g_{\underline {\sigma} \sigma x} (\eta) -\Psi (\bar {\eta}_{\sigma x_0,l})\right |
\\
& \le g^* \frac {1} { | \{ \underline { \sigma } | |\underline {\sigma }| \le l \} |}   \sum_{ | \underline {\sigma } | \le l} \eta_{\underline {\sigma} \sigma x} (\eta) +g^* \bar {\eta}_{\sigma x_0,l}
\\
&\le (g^* |D_{x_0}| +g^*) \bar {\eta}_{\sigma x_0,l}.
\end{align*}

\noindent Together with Lemma 4.4, to prove the one block estimate, it suffices to show that, for every $C_1>0$, 

\begin{equation} \label{cutoff}
\varlimsup_{l \to \infty} \varlimsup_{N \to \infty} \sup_{D_N (f) \le \frac{C |V_N| } {N^2}} \int \frac{1}{ |\Gamma_N| } \sum_{\sigma \in \Gamma_N} V_{\sigma x,l} (\eta) 1_{ \{ \bar{\eta} _{\sigma x_0,l} \le C_1\} } f(\eta) \nu_{\alpha^*} (d\eta). 
\end{equation}

\noindent Since $\nu_{\alpha^*}$ is $\Gamma_N$ invariant, we have the integral part above is equal to

$$\int V_{x,l} (\eta) 1_{ \{ \bar{\eta}_{x_0,l} \le C_1 \}} \bar{f} (\eta) \nu_{\alpha^*}$$

\noindent where $\bar{f} (\eta) := \frac{1}{ |\Gamma_N| } \sum_{\sigma \in \Gamma_N} \sigma f(\eta)$. Notice that $V_{x,l} (\eta) 1_{ \{ \bar{\eta}_{x_0,l} \le C_1 \}}$ only depends on the coordinates $\{ \eta_x,\ x\in B(D_{x_0}, l)\}$. For fixed $l$, we will consider a subgraph $\Lambda_l:=(V_l, E_l)$, where $V_l:=B(D_{x_0}, l)$ and $E_l:= \{ e \in E | oe,te \in V_l\}$. Let $X^l:=\mathbb{N}^{V_l}$ be the configuration space and $\nu_{\alpha^*}^l$ be the restriction of $\nu_{\alpha*}$ to $X^l$, i.e., $\nu_{\alpha^*}^l (\xi):=\nu_{\alpha^*} \{ \eta | \eta_x=\xi_x, x \in V_l \}$. For a density function $f$, we represent by $f_l$ the conditional expectation of $f$ to the $\sigma$-algebra generated by $\{\eta_x, x \in V_l\}$, i.e.,

$$f_l(\xi):=\frac{1}{\nu_{\alpha^*}^l (\xi)} E_{\nu_{\alpha^*}} [f(\eta) 1_{ \{ \eta| \eta_x=\xi_x,x \in V_l \} }],\ \ \ \ \ \ \ \ \ \xi \in X^l.$$

\noindent With these notations, we can rewrite (\ref{cutoff}) as

\begin{equation} \label{cutoff2}
\varlimsup_{l \to \infty } \varlimsup_{N \to \infty} \sup_{D_N(f) \le \frac {C |V_N| } {N^2}} \int V_{x,l} (\xi) 1_{ \{ \bar{\xi} _{x_0,l} \le C_1 \}} \bar{f}_l (\xi) \nu_{\alpha^*}^l (d \xi).
\end{equation}

\noindent Next we estimate the Dirichlet form of $\bar{f}_l$. For each $e \in E$, define 

\begin{equation}
L_{oe,te}f(\eta) := \frac{1}{2} \{ g(\eta_{oe}) [f(\eta^e) -f(\eta) ] +g(\eta_{te}) [f(\eta^{\bar{e}}) -f(\eta)]\},
\end{equation}

\begin{equation}
I_{oe,te} (f):= -\langle L_{oe,te} \sqrt{f}, \sqrt{f} \rangle _{\nu_{\alpha^*}}=\frac{1}{2} \int g(\eta _{oe}) \{ \sqrt{f(\eta^e)} -\sqrt{f (\eta)} \}^2 \nu_{\alpha^*} (d\eta).
\end{equation}
\noindent Then we have 

$$L_Nf(\eta)=\sum_{ e \in E_N} p(e)L_{oe,te} f(\eta),$$

$$ D_N(f)=\sum_{ e \in E_N} p(e)I_{oe,te} (f).$$

\noindent Let us restrict the above definition to $X^l$: for every density function $h: X^l \to \mathbb{R}$, 

 $$I_{oe,te}^l( h):=\frac{1}{2} \int g(\xi _{oe}) \{ \sqrt{h(\xi^e)} -\sqrt{h (\xi)} \}^2 \nu_{\alpha^*}^l (d\xi),$$
 
 $$D^l(h):=\sum_{e \in E_l} I_{oe,te}^l (h).$$
 By Cauchy-Schwarz inequality,
 
 $$I_{oe,te}^l (\bar{f} _l) \le I_{oe,te} (\bar{f}).$$
Note that $I_{oe,te} (\sigma f) =I_{\sigma oe, \sigma te} (f)$, which implies $D_N( \sigma f) =D_N(f)$ and combine with the convexity of the Dirichlet form,  
 it holds that,
 
 $$D^l( \bar{f}_l) \le \sum_{e \in E_l} I_{oe,te} (\bar{f}) \le | \{ \sigma | |\sigma| \le l \}| \frac{1}{ p^*|\Gamma_N| } D_N(\bar{f}) \le \frac{C |D_{x_0}| (2l+1)^d } {p^*N^2 } $$
 where $p^*:=\min_{e \in E} p(e)>0$ is a constant. Since we cut off the density by the indicator function, we can restrict the supremum to the densities concentrated on the set $\{ \xi | \bar{\xi} _{x_0, l} \le C_1 \}$. This subset of $\mathcal{P}{(X^l)}$ is compact for the weak topology. Furthermore, by the continuity of the Dirichlet form, as $N \to \infty$, (\ref{cutoff2}) is bounded above by 

 \begin{equation}
 \varlimsup_{ l \to \infty} \sup_{D^l (f) =0} \int V_{x,l} (\xi) 1_{ \{ \bar{\xi}_{x_0,l} \le C_1 \} } f(\xi) \nu^l_{\alpha^*} (\xi) =0.
 \end{equation}
 
 \noindent For each $j\ge 0$, let $\nu^{l,j} $ be the measure $\nu_{\alpha^*}^l$ conditioned to the hyperplane $\{ \xi | \sum_{x \in V_l} \xi_x =j \}$, i.e.,
 
 $$\nu^{l,j} (\cdot):=\frac{\nu^l_{\alpha^*} (\cdot)  } {\nu^l_{\alpha^*} \{\xi | \sum_{x \in V_l} \xi_x =j \}}.$$

\noindent Note that $D^l (f) =0$ implies that $f$ is constant on each hyperplane $\{\xi | \sum_{x \in V_l} \xi_x =j \}$ for each $j \le 0$.

\noindent We have that

\begin{align*}
\int V_{x,l} (\xi) 1_{ \{ \bar{\xi}_{x_0,l} \le C_1 \} } f(\xi) \nu^l_{\alpha^*} (\xi) &= \sum_{j=0}^{C_1 |V_l |} \int_{\{\sum_{x \in V_l} \xi_x =j \}} V_{x,l} (\xi)  f(\xi) \nu^l_{\alpha^*} (\xi)
\\
&= \sum_{j=0}^{C_1 |V_l |} \int_{\{\sum_{x \in V_l} \xi_x =j \}}  f(\xi) \nu^l_{\alpha^*} (\xi) \int V_{x,l}(\xi) \nu^{l.j} (d \xi).
\end{align*}

\noindent Since $\sum_{j=0}^{C_1 |V_l |} \int_{\{\sum_{x \in V_l} \xi_x =j \}}  f(\xi) \nu^l_{\alpha^*} (\xi)=1$, it is enough to show that

\begin{equation}
\varlimsup_{l \to \infty} \sup_{j \le C_1 |V_l |} \int V_{x,l}(\xi) \nu^{l.j} (d \xi) =0.
\end{equation}

 For a fixed positive integer $k$, define $A:= A_{l-k} \cap \{ (2k+1) \sigma, \sigma \in \Gamma \}:=\{ \sigma_i\}_{i=1}^n$, where $A_m:= \{ | \sigma | \le m \}$ for each positive integer $m$. For each $1\le i \le n$, let $B_i=\sigma_i +A_k$. Then $\{B_i\}_{i=1}^n$ are pairwise disjoint and belong to $A_l$. Let $B_0:=A_l \setminus \cup_{i=1}^n B_i $, then $|B_0| \le |A_l|-|A_{l-k}| \le c(k){l^{d-1}}$. Notice that $\nu^{l,j}$ is concentrated on the configurations with $j$ particles, then 

\begin{equation*}
\begin{aligned}
\int V_{x,l} (\xi) \nu^{l,j} (d \xi) &= \int \left| \frac {1} { | \{\sigma | | \sigma | \le l \} | } \sum_{ | \sigma | \le l} g(\eta_{\sigma x}) -\Psi (\bar {\xi} _{x_0,l} ) \right| \nu^{l,j} (d \xi) 
\\
&\le \sum _{i=0}^n \frac { |B_i| } { | A_l| } \int \left |\frac {1} { | B_i | } \sum_{\sigma \in B_i} g(\xi_{\sigma x}) -E_{\nu_{\frac{j}{ | V_l | }}} [g(\xi _{x_0})] \right| \nu^{l,j} (d \xi)
\\
&=\sum _{i=1}^n \frac { |B_i| } { | A_l| } \int \left |\frac {1} { | B_i | } \sum_{\sigma \in B_i} g(\xi_{\sigma x}) -E_{\nu_{\frac{j}{ | V_l | }}} [g(\xi _{x_0})] \right| \nu^{l,j} (d \xi)
\\
&\ \ \ \ \ \ + \frac { |B_0| } { | A_l| } \int \left |\frac {1} { | B_0 | } \sum_{\sigma \in B_0} g(\xi_{\sigma x}) -E_{\nu_{\frac{j}{ | V_l | }}} [g(\xi _{x_0})] \right| \nu^{l,j} (d \xi).
\end{aligned}
\end{equation*}

\noindent The last term tends to 0 as $l \to \infty$ since $|B_0| \le c(k){l^{d-1}}$. Since $\{\xi_x, x \in B_i \}_{i=1}^n$ have the same distribution, the previous summation is bounded by

$$ \int \left | \frac {1} { | \{\sigma | | \sigma | \le k \} | } \sum_{ | \sigma | \le k} g(\eta_{\sigma x}) -E_{\nu_{\frac{j}{|V_l|}}} [g(\xi_{x_0})] \right| \nu^{l,j} (d \xi) .$$

\noindent By the equivalence of ensembles(see Appendix 2 of \cite{KL1999}), as $l \to \infty$ and $\frac{j}{|V_l|} \to \alpha$, the above integral converges (uniformly in $\alpha$ on each interval of $\mathbb{R}_+$) to 

$$ \int \left | \frac {1} { | \{\sigma | | \sigma | \le k \} | } \sum_{ | \sigma | \le k} g(\eta_{\sigma x}) -E_{\nu_{\alpha}} [g(\xi_{x_0})] \right| \nu_{\alpha} (d \xi). $$

\noindent By the law of large numbers, the last integral converges (uniformly in $\alpha$ on each interval of $\mathbb{R}_+$) to 0 as $k \to \infty$. It completes the proof of the one block estimate Lemma 4.2.

\subsection{Proof of two blocks estimate lemma} 

As in Section 4.3, for fixed $l$, we consider the subgraph $\Lambda_l:=(V_l, E_l)$, where $V_l:=B(D_{x_0}, l)$ and $E_l:= \{ e \in E | oe,te \in V_l\}$. Let $X^l \times X^l:=\mathbb{N}^{V_l} \times \mathbb{N}^{V_l}$ be the configuration space and $\nu_{\alpha^*}^l \otimes \nu_{\alpha^*}^l  $ be the product measure on $X^l \times X^l$. For a density function $f : Z_N \to \mathbb{R}_+$, $f_{\sigma ,l}$ is the conditional expectation of $f$ with respect to the $\sigma $-algebra generated by $\{\eta_x, x \in V_l \cup \sigma V_l\}$, i.e., for every $(\xi,\xi^\prime) \in X^l \times X^l$,

$$ f_{\sigma,l} (\xi, \xi^\prime)=\frac{1}{ \nu_{\alpha^*}^l \otimes  \nu_{\alpha^*}^l (\xi,\xi^\prime)} E_{\nu_{\alpha^*}} [ f (\eta) 1_{ \{ \eta\ |\ \eta_x=\xi_x, \eta_{\sigma x} =\xi^\prime_x,\ x \in V_l \} }].$$

\noindent For a density $h: X^l \times X^l \to \mathbb{R}_+$, define the Dirichlet form of $h$ by

\begin{equation}
D^{l,l} (h):=\sum_{e \in E_l} (I_{oe,te}^{1,l} +I_{oe,te}^{2,l}) +I_{x_0,x_0^\prime}^l(h),
\end{equation}

\noindent where 
\begin{equation}
\begin{aligned}
&I_{oe,te}^{1,l}:=\frac{1}{2} \int g(\xi_{oe}) [\sqrt {h(\xi^e,\xi^\prime)} -\sqrt {h(\xi,\xi^{\prime })} ]^2 \nu_{\alpha^*}^l \otimes  \nu_{\alpha^*}^l (d\xi,d\xi^\prime),
\\
&I_{oe,te}^{2,l}:=\frac{1}{2} \int g(\xi_{oe}) [\sqrt {h(\xi,\xi^{\prime e})} -\sqrt {h(\xi,\xi^{\prime })} ]^2 \nu_{\alpha^*}^l \otimes  \nu_{\alpha^*}^l (d\xi,d\xi^\prime),
\\
&I_{x_0,x_0^\prime}^ l (h):=\frac{1}{2 } \int g(\xi_{x_0}) [ \sqrt {h(\xi ^{x_0,-} ,\xi ^{\prime x_0,+})} -\sqrt{h(\xi, \xi^\prime)}] ^2  \nu_{\alpha^*}^l \otimes  \nu_{\alpha^*}^l (d\xi,d\xi^\prime),
\end{aligned}
\end{equation}
where 
\begin{equation*}
\xi^{x,\pm}_z=
\left\{
\begin{aligned}
& \xi_{x} \pm 1 \ \ \ \ z=x
\\
& \xi_z\ \ \ \ \ otherwise.
\end{aligned}
\right.
\end{equation*}

\noindent For any $x,y \in V_N$, define 

\begin{equation}
I_{x,y}:=\frac{1}{2} \int g(\eta_x)  \{ \sqrt{f(\eta^{x,y})} -\sqrt{f (\eta)} \}^2 \nu_{\alpha^*} (d\eta),
\end{equation}

\noindent where $\eta^{x,y}$ stands for the configuration of $\eta$ where a particle jumped from $x$ to $y$, i.e.,

\begin{equation*}
\eta^{x,y}_z=
\left\{
\begin{aligned}
& \eta_x -1 \ \ \ \ z=x
\\
&\eta_y+1\ \ \ \ z=y
\\
& \eta_z\ \ \ \ \ otherwise.
\end{aligned}
\right.
\end{equation*}
\noindent We first give an lemma needed later.

\begin{lemma}
There exists a constant $c_1>0$ such that for all $\sigma \in \Gamma$, it holds that
\begin{equation}
d(x_0,\sigma x_0) \le c_1 |\sigma|,
\end{equation}
where $d$ is the graph distance of $X$.
\end{lemma}
\noindent This lemma is easy to show by taking $c_1=\sup_{|\sigma|=1} d(x_0,\sigma x_0)$ and induction, so we omit the proof.

\begin{lemma} \label{XY}
For every $\sigma \in \Gamma_N$, 
\begin{equation}
I_{x_0, \sigma x_0} (\bar{f}) \le d(x_0, \sigma x_0) ^2 \frac{D_N(\bar{f})} { p^*| \Gamma_N | }.
\end{equation}
\end{lemma}

\noindent Proof: For $x_0, \sigma x_0$, there exists a path $c=(e_1,\dots,e_n)$ such that $x_0=oe_1, \dots, te_n=\sigma x_0$ and $n=d(x_0,\sigma x_0)$. After a change of variable $\eta=\xi^{x_0,+}$, we have

\begin{align*}
 I_{x_0, \sigma x_0} (\bar{f}) &=\frac{1}{2} \Psi(\alpha^*) \int \left ( \sqrt{ \bar{f} ( \eta^ { \sigma x_0,+})} -\sqrt { \bar{ f } ( \eta^{x_0,+})} \right)^2 \nu_{\alpha^*} (d\eta)
 \\
 &=\frac{1}{2} \Psi(\alpha^*) \int \left ( \sum_{k=1}^n (\sqrt{ \bar{f} ( \eta^ { te_k,+})} -\sqrt { \bar{ f } ( \eta^{oe_k,+})}) \right )^2\nu_{\alpha^*} (d\eta)
 \\
 &\le \frac{1}{2} \Psi(\alpha^*) n \sum_{k=1}^n  \int \left ( \sqrt{ \bar{f} ( \eta^ { te_k,+})} -\sqrt { \bar{ f } ( \eta^{oe_k,+})}) \right )^2\nu_{\alpha^*} (d\eta)
 \\
 &= n\sum_{k=1}^n I_{oe_k,te_k} (\bar{f}) \le d(x_0, \sigma x_0)^2 \frac{D_N(\bar{f})} {p^*|\Gamma_N|}.
 \end{align*}

\noindent Proof of Lemma \ref{two}:  As in the proof of the one block estimate, we can cut off of large density first. It is enough to show that 

\begin{equation} \label{cut}
\begin{aligned}
\varlimsup_{l \to \infty} \varlimsup_{\epsilon \to \infty} \varlimsup_{N \to \infty} &\sup_{D_N (f) \le \frac{C | V_N |}{N^2}} \sup_{2l \le |\sigma| \le \epsilon N} \\
&\int | \bar{\eta}_{x_0,l} -\bar{\eta}_{\sigma x_0,l}| 1_{ \{\bar{\eta}_{x_0,l} +\bar{\eta}_{\sigma x_0,l} \le A \} } \bar{f}(\eta) \nu_{\alpha^*} (d\eta) =0.
\end{aligned}
\end{equation}

\noindent (\ref{cut}) can be written as 

\begin{equation}
\begin{aligned}
\varlimsup_{l \to \infty} \varlimsup_{\epsilon \to \infty}& \varlimsup_{N \to \infty}\sup_{D_N (f) \le \frac{C | V_N |}{N^2}} \sup_{2l \le |\sigma| \le \epsilon N} \\
& \int  | \bar{\xi}_{x_0,l} -\bar{\xi^\prime}_{ x_0,l}| 1_{ \{\bar{\xi}_{x_0,l} +\bar{\xi^\prime}_{x_0,l} \le A \} } \bar{f}_{\sigma, l} (\xi,\xi^\prime) \nu_{\alpha^*}^l \otimes  \nu_{\alpha^*}^l  (d\xi, d\xi^\prime) =0.
\end{aligned}
\end{equation}
As in the proof of Lemma \ref{one}, the next step consists in the estimation of the Dirichlet form $D^{l,l}(\bar{f}_{\sigma,l})$. 
Note that
\begin{equation}
\begin{aligned}
\sum_{e \in E_l} (I_{oe,te}^{1,l} (\bar{f}_{\sigma,l}) +I_{oe,te}^{2,l} (\bar{f}_{\sigma,l})) \le \sum_{e \in E_l \cup \sigma E_l} I_{oe,te} (\bar{f}) \le 2 |E_l| \frac{D_N(\bar{f})}{ p^*|\Gamma_N|} \le \frac{2 C |E_l| |V_0|}{p^*N^2}
\end{aligned}
\end{equation}
and by Lemma \ref{XY}

\begin{equation}
\begin{aligned}
I_{x_0,x_0^\prime}^l ( \bar{f}_{\sigma,l}) \le I_{x_0,\sigma x_0} (\bar{f}) \le d(x_0, \sigma x_0) ^2 \frac{D_N( \bar{ f} )} {p^* | \Gamma_N | } \le p^{*-1}c_1^2 \epsilon^2 C |V_0| .
\end{aligned}
\end{equation}
For the same reason of Lemma \ref{one}, it is enough to prove that

\begin{equation}
\varlimsup_{l \to \infty} \sup_{D^{l,l}(f) =0} \int  | \bar{\xi}_{x_0,l} -\bar{\xi^\prime}_{ x_0,l}| 1_{ \{\bar{\xi}_{x_0,l} +\bar{\xi^\prime}_{x_0,l} \le A \} } f (\xi,\xi^\prime) \nu_{\alpha^*}^l \otimes  \nu_{\alpha^*}^l  (d\xi, d\xi^\prime) =0.
\end{equation}
The proof will be completed in the way as mentioned in the proof of Lemma \ref{one}.

\section{ standard realization }
In this section, we apply our results to find the standard realization of the crystal lattices. For a crystal lattice $X=( V , E )$, fix a harmonic realization $\Phi_0$ with lattice group $\phi_0(\Gamma)=\{\sum_{i=1}^d k_iu_i\ | k_i \in \mathbb{Z}\ \}$. Since the diffusion matrix $\mathbb{D}_{\Phi_0}$ is strictly positive definite, all eigenvalues $\{\lambda_1,\cdots,\lambda_d\}$ of $\mathbb{D}_{\Phi_0}$ are  strictly positive. 
We can write $\mathbb{D}_{\Phi_0}$ as

\begin{equation}
\mathbb{D}_{\Phi_0}=P^T diag(\lambda_1,\dots,\lambda_d) P.
\end{equation}
Here $P$ is an orthogonal matrix and $diag(\lambda_1,\dots,\lambda_d)$ is a diagonal matrix. 
 
The following proposition tells us how to get the standard realization from a fixed harmonic realization.

\begin{proposition} \label{standard}
For all $\Phi^\prime$ with $vol(D_{\Phi^\prime})=vol(D_\Phi)$, it holds that
 $$E(\Phi)\le E(\Phi^\prime)$$
 where $\Phi$ is the harmonic realization with the lattice group $\{\sum_{i=1}^d k_iA u_i\ | k_i \in \mathbb{Z}\ \}$ and $A$ is given by
 \begin{equation}
 A= diag\left ( \left (\frac{\lambda_1\cdots\lambda_d}{\lambda_1^d} \right )^{\frac{1}{2d}} , \dots, 
\left (\frac{\lambda_1\cdots\lambda_d}{\lambda_d^d} \right )^{\frac{1}{2d}} \right) P.
\end{equation}
\end{proposition}
\noindent Proof:
For any harmonic realization $\Phi$ with lattice group $\phi(\Gamma)=\{ \sum _{i=1}^d k_i\widetilde{u}_i\ |\ k_i \in \mathbb{Z}\}$, let $A$ be the basis transformation from $\{u_1,\dots,u_d\} $ to $\{\widetilde{u}_1,\dots,\widetilde{u}_d\}$. By Theorem \ref{exclusion} and Proposition \ref{be}, we have that
\begin{equation*}
\mathbb{D}_\Phi=A \mathbb{D}_{\Phi_0}A^T.
\end{equation*}
Note that the energy is nothing but the trace of diffusion matrix, i.e.,

\begin{equation*}
E(\Phi)=tr(\mathbb{D}_\Phi)=tr(A \mathbb{D}_{\Phi_0}A^T).
\end{equation*}
Thus, fix the volume of fundamental parallelotope, to find the standard realization, it suffices to find $A$ with $|A|=1$ such that
\begin{equation*}
tr(A \mathbb{D}_{\Phi_0}A^T)=\min_{|A|=1}tr(A \mathbb{D}_{\Phi_0}A^T).
\end{equation*}
By Schwarz inequality, we have that

\begin{equation*}
tr(A \mathbb{D}_{\Phi_0}A^T)\ge d \sqrt[d]{|A \mathbb{D}_{\Phi_0}A^T|}=d \sqrt[d]{| \mathbb{D}_{\Phi_0}|}.
\end{equation*}
Then the above inequality takes the equal when $A$ equal to 

\begin{equation*}
A= diag\left ( \left (\frac{\lambda_1\cdots\lambda_d}{\lambda_1^d} \right )^{\frac{1}{2d}} , \dots, 
\left (\frac{\lambda_1\cdots\lambda_d}{\lambda_d^d} \right )^{\frac{1}{2d}} \right) P.
\end{equation*}
As mentioned in Remark \ref{har}, the harmonic realization with lattice group $\{\sum_{i=1}^d k_iAu_i\ | k_i \in \mathbb{Z}\ \}$ is also the standard realization. The proof is completed.

\section{examples}
In this section, we give two concrete examples where we can apply our main results. In particular, both models are given as a process on the square lattice with inhomogeneous jump rates. We study theses models as a homogeneous model on a crystal lattice and obtain explicit hydrodynamic equations.

\noindent
\textbf{Example 1.} Consider the exclusion process on discrete torus $\mathbb{T}_{2N}:=\{ 0,1,\dots, 2N-1\}$ with generator:

\begin{equation}
L_Nf(\eta):=\sum_{x \in \mathbb{ T }_{ 2N } } p(x,x+1) \left\{ f(\eta^{x,x+1}) -f(\eta) \right \},
\end{equation}
where 

\begin{equation}
p(x,x+1)=\left\{
\begin{aligned}
&\alpha\ \ \ \ \  x\  \ even
\\
\\
&\beta\ \ \ \ \ x\  \ odd.
\end{aligned}
\right.
\end{equation}

\noindent We can regard it as a process on the crystal $X=(V,E)$, where $V=\mathbb{Z}, E=\{ (x,x+1),(x+1,x)\ |\ z \in \mathbb{Z} \}$ and the group $\Gamma=\mathbb{Z}$ with the group action: $\sigma x:=2\sigma +x$ and for $\sigma \in \Gamma$. Then $V_0=\{0,1\}$, $\mathbb{T}_\phi=[0,2)$, $\Gamma_N=\{0,1,\cdots,N-1\}$, $V_N=\{0,1,\dots,2N-1\}$ with $2N=0$ and $E_N=\{(x,x+1),(x+1,x)\ |\ x \in V_N\}$. The generator is given by

\begin{equation}
L_N f ( \eta ):=\sum_{ e \in E_N } p( e ) \{  f ( \eta^e ) - f ( \eta ) \}.
\end{equation}

\noindent Take the realization $\Phi$ as in example 1b, the associated empirical process is

\begin{equation*}
 \pi_t^{\Phi,N} ( du )=\frac{1}{2N} \sum_{x \in \{0,1\}} \sum_{y=1}^{N-1} \eta_{2y+x}(t) \delta_{2y+x/N} (du).
\end{equation*}
Note that this realization is not harmonic for the lattice group $\{\phi(\sigma)=2\sigma,\ \sigma \in \Gamma\}$ and the weight function $p(\cdot)$. We obtain the harmonic realization $\Phi_h$(See Figure 5) associated to the lattice group $\phi(\Gamma)$  by (\ref{harmonic}). Precisely, 
\begin{equation}
\Phi_h(\sigma 0):=0+\phi(\sigma),\  {\Phi_h(\sigma 1):=\frac{2\beta}{\alpha+\beta}}+\phi(\sigma).
\end{equation}

\noindent By Theorem \ref {exclusion}, we have that 
\begin{equation}
\varlimsup_{N \to \infty} P_{\mu^N}\left[\left| \frac{1}{2N}  \sum_{ x \in \mathbb{T}_{2N}}G \left ( \frac{x}{N} \right ) \eta_x(t)-\frac{1}{2} \int_0^2 G(u)\rho(t,u)du \right|>\delta \right]=0
\end{equation}
for every $\delta>0$, every continuous $G: [0,2) \to \mathbb{R} $ and $\rho(t,u)$ is the unique weak solution of
\begin{equation}
 \left\{
 \begin{aligned}
 & \partial_t\rho = \frac{4 \alpha \beta}{ \alpha+\beta} \partial^2_ x  \rho 
  \\
 & \rho(0,\cdot)=\rho_0(\cdot)
  \end{aligned}
  \right.
   \end{equation}
where $ \frac{4\alpha \beta }{ \alpha+\beta}$ is just the diffusion coefficient (matrix) of $ \Phi_h$.
\begin{remark}
Note that if $\beta$ is much bigger than $\alpha$, the diffusion coefficient is close to $4\alpha$. This means: from microscopic view, particle jumps very fast between site 1 and site 2; in macroscopic view, the diffusion speed is close to $4\alpha$, where $\alpha$ is the jump rate between site 0 and site 1.  
\end{remark}

 \begin{figure} \label{Ex1}
      	\begin{tikzpicture} [xscale = 1, yscale = 1]
	
	\draw [ black, very thick]  (-1, -3) -- (8,-3) ;  	   
      	\draw [very thick, fill = white] (2, -3) circle[radius = 2.5pt] node[below = 0.2cm] {$\Phi ( 1 )$};
      	\draw [very thick, fill = white] (5, -3) circle[radius = 2.5pt] node[below = 0.2cm] {$\Phi ( 3 )$};
      	\draw[very thick,fill=black]  ( 0.5 , -3) circle[radius = 2.5pt] node[below = 0.2cm] {$\Phi(0)$};
      	\draw[very thick, fill=black] (3.5, -3) circle[radius = 2.5pt] node[below = 0.2cm] {$\Phi(2)$}; 
      	\draw[fill=black,very thick] (6.5, -3) circle[radius = 2.5pt] node[below = 0.2cm] {$\Phi(4)$};
      	\draw [  very thick, black] (0.55 , -2.8) arc [start angle = 150, end angle = 30, radius = 0.8];
      	\draw [  very thick, black] (2.05 , -2.8) arc [start angle = 150, end angle = 30, radius = 0.8];
      	\draw [  very thick, black] (3.55 , -2.8) arc [start angle = 150, end angle = 30, radius = 0.8];
      	\draw [  very thick, black] (5.05 , -2.8) arc [start angle = 150, end angle = 30, radius = 0.8];
      	\draw (1.25, -2.35) node[above] { \small $ 1 $ };
      	\draw (4.25, -2.35) node[above] { \small $ 1 $ }; 
      	\draw (2.75, -2.35) node[above] { \small $ 1 $ };
      	\draw (5.75, -2.35) node[above] { \small $ 1 $ };
	
	\draw (1.25, -3) node[above] { \small $ \alpha $ };
	\draw (2.75, -3.05) node[above] { \small $ \beta $ };
      	\draw (4.25, -3) node[above] { \small $ \alpha $ }; 
      	\draw (5.75, -3.05) node[above] { \small $ \beta $ };

	\draw [ black, very thick]  (-1, -6) -- (8,-6) ;  
	\draw[very thick,fill=black]  ( 0.5 , -6) circle[radius = 2.5pt] node[below = 0.2cm] {$\Phi_h(0)$};
	\draw[very thick,fill=black]  ( 3.5 , -6) circle[radius = 2.5pt] node[below = 0.2cm] {$\Phi_h(2)$};
	\draw[very thick,fill=black]  ( 6.5 , -6) circle[radius = 2.5pt] node[below = 0.2cm] {$\Phi_h(4)$};
	\draw [very thick, fill = white] (2.3, -6) circle[radius = 2.5pt] node[below = 0.2cm] {$\Phi_h(1)$};
	\draw [very thick, fill = white] (5.3, -6) circle[radius = 2.5pt] node[below = 0.2cm] {$\Phi_h(3)$};
	
	\draw [  very thick, black] (0.55 , -5.8) arc [start angle = 135, end angle = 45, radius = 1.15];
      	\draw [  very thick, black] (2.35 , -5.8) arc [start angle = 158, end angle = 22, radius = 0.55];
      	\draw [  very thick, black] (3.55 , -5.8) arc [start angle = 135, end angle = 45, radius = 1.15];
      	\draw [  very thick, black] (5.35 , -5.8) arc [start angle = 158, end angle = 22, radius = 0.55];

        \draw (1.35, -5.2) node[above] { \small $ \frac{2 \beta }{\alpha+\beta} $ };
      	\draw (4.35, -5.2) node[above] { \small $ \frac{ 2 \beta }{ \alpha+\beta} $ }; 
      	\draw (2.85, -5.2) node[above] { \small $ \frac{2 \alpha }{\alpha+\beta} $ };
      	\draw (5.85, -5.2) node[above] { \small $\frac{ 2 \alpha }{ \alpha+\beta} $ };
	
	\draw (1.4, -6) node[above] { \small $ \alpha $ };
	\draw (2.9, -6.05) node[above] { \small $ \beta $ };
      	\draw (4.4, -6) node[above] { \small $ \alpha $ }; 
      	\draw (5.9, -6.05) node[above] { \small $ \beta $ };

      	\end{tikzpicture}
	\caption{The images of $\Phi$(above) and $\Phi_h$(below). The numbers in the first line show the length between vertices in each realization, while those in the second line show the jump rates between vertices.}
      \end{figure}
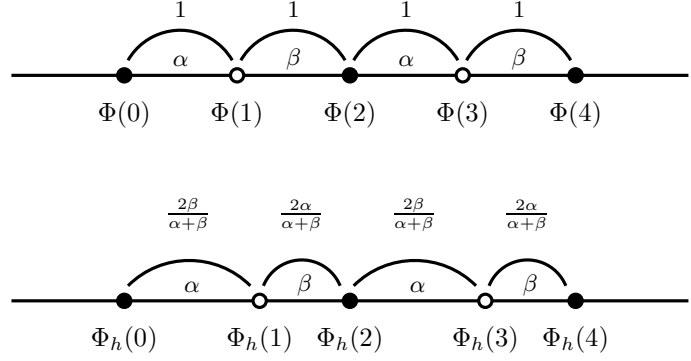

\noindent\textbf{Example 2.} Let $ \Gamma:= \{ k_1 (2,0 )+ k_2 (0,1) | k_1, k_2 \in \mathbb{Z}\}$ be a lattice group on $\mathbb{Z}^2$. Let $ p( \cdot, \cdot ): \mathbb{Z}^2 \times \mathbb{ Z}^2 \to [0,1]$ be a weight function satisfying:

(i) $p(x,y)=p(y,x)$ for all $ x, y \in \mathbb{Z}^2$.

(ii) $p( \sigma x ,\sigma y)=p(x,y)$ for all $x,y \in \mathbb{Z}^2$ and all $\sigma \in \Gamma$, where $ \sigma x:=x+\sigma$.

(iii) $p(x,y)=0 $ if $|x-y |:=|x_1-y_1|+| x_2 -y_2 | \neq 1$.

(iv) $p((0,0), (1,0 ) )=\frac{1}{6},p( ( 0, 0 ) ,( 0,1 ) )=\frac{1}{3}, p(( 0, 0 ), ( -1,0 ) ) =\frac{1}{2}, p( (0,0) , ( 0,-1))=0$.

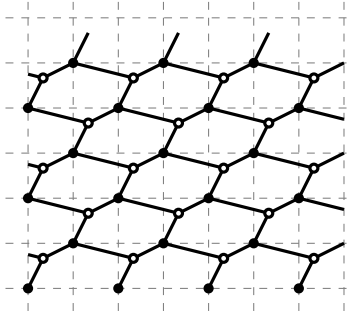
\begin{figure} \label{2}
      \begin{tikzpicture} [xscale=0.6,yscale=0.6]

      \draw [dashed, gray] (-0.5, 0 ) -- (7.5, 0);
       \draw [dashed, gray] (-0.5, 1 ) -- (7.5, 1);
        \draw [dashed, gray] (-0.5, 2 ) -- (7.5, 2);
         \draw [dashed, gray] (-0.5, 3 ) -- (7.5, 3);
          \draw [dashed, gray] (-0.5, 4 ) -- (7.5, 4);
           \draw [dashed, gray] (-0.5, 5 ) -- (7.5, 5);
            \draw [dashed, gray] (-0.5, 6 ) -- (7.5, 6);
             
              \draw [dashed, gray] (0,-0.5 ) -- ( 0,6.5);
              \draw [dashed, gray] (1,-0.5 ) -- ( 1,6.5);
              \draw [dashed, gray] (2,-0.5 ) -- ( 2,6.5);
              \draw [dashed, gray] (3,-0.5 ) -- ( 3,6.5);
              \draw [dashed, gray] (4,-0.5 ) -- ( 4,6.5);
              \draw [dashed, gray] (5,-0.5 ) -- ( 5,6.5);
              \draw [dashed, gray] (6,-0.5 ) -- ( 6,6.5);
               \draw [dashed, gray] (7,-0.5 ) -- ( 7,6.5);

                \draw[black,very thick] (0,0)--(1/3,2/3);
            \draw[black,very thick] (2,0)--(7/3,2/3);
             \draw[black,very thick] (4,0)--(13/3,2/3);
              \draw[black,very thick] (6,0)--(19/3,2/3);
              
                \draw[black,very thick] (1,1)--(1/3,2/3);
                 \draw[black,very thick] (3,1)--(7/3,2/3);
             \draw[black,very thick] (5,1)--(13/3,2/3);
               \draw[black,very thick] (7,1)--(19/3,2/3);
                             
             \draw[black,very thick] (0,2)--(1/3,8/3);
            \draw[black,very thick] (2,2)--(7/3,8/3);
             \draw[black,very thick] (4,2)--(13/3,8/3);
              \draw[black,very thick] (6,2)--(19/3,8/3);
              
                \draw[black,very thick] (1,3)--(1/3,8/3);
                 \draw[black,very thick] (3,3)--(7/3,8/3);
             \draw[black,very thick] (5,3)--(13/3,8/3);
               \draw[black,very thick] (7,3)--(19/3,8/3);
               
                 \draw[black,very thick] (1,1)--(7/3,2/3);
                 \draw[black,very thick] (3,1)--(13/3,2/3);
             \draw[black,very thick] (5,1)--(19/3,2/3);
             
               \draw[black,very thick] (1,3)--(7/3,8/3);
                 \draw[black,very thick] (3,3)--(13/3,8/3);
             \draw[black,very thick] (5,3)--(19/3,8/3);
             
              \draw[black,very thick] (1,1)--(4/3,5/3);
               \draw[black,very thick] (3,1)--(10/3,5/3);
                \draw[black,very thick] (5,1)--(16/3,5/3);
                
                \draw[black,very thick] (1/3,2/3)--(0,3/4);
                  \draw[black,very thick] (2,2)--(4/3,5/3);
               \draw[black,very thick] (4,2)--(10/3,5/3);
                \draw[black,very thick] (6,2)--(16/3,5/3);
                
                  \draw[black,very thick] (0,2)--(4/3,5/3);
               \draw[black,very thick] (2,2)--(10/3,5/3);
                \draw[black,very thick] (4,2)--(16/3,5/3);
                
                 \draw[black,very thick] (0,4)--(1/3,14/3);
            \draw[black,very thick] (2,4)--(7/3,14/3);
             \draw[black,very thick] (4,4)--(13/3,14/3);
              \draw[black,very thick] (6,4)--(19/3,14/3);
              
              \draw[black,very thick] (1,5)--(1/3,14/3);
                 \draw[black,very thick] (3,5)--(7/3,14/3);
             \draw[black,very thick] (5,5)--(13/3,14/3);
               \draw[black,very thick] (7,5)--(19/3,14/3);
               
                \draw[black,very thick] (1/3,8/3)--(0,11/4);
                \draw[black,very thick] (1/3,14/3)--(0,19/4);
                
                 \draw[black,very thick] (1,3)--(4/3,11/3);
               \draw[black,very thick] (3,3)--(10/3,11/3);
                \draw[black,very thick] (5,3)--(16/3,11/3);
                
                   \draw[black,very thick] (2,4)--(4/3,11/3);
               \draw[black,very thick] (4,4)--(10/3,11/3);
                \draw[black,very thick] (6,4)--(16/3,11/3);
                
                 \draw[black,very thick] (1,5)--(4/3,17/3);
               \draw[black,very thick] (3,5)--(10/3,17/3);
                \draw[black,very thick] (5,5)--(16/3,17/3);
                
                   \draw[black,very thick] (0,4)--(4/3,11/3);
               \draw[black,very thick] (2,4)--(10/3,11/3);
                \draw[black,very thick] (4,4)--(16/3,11/3);
                
                  \draw[black,very thick] (1,5)--(7/3,14/3);
                 \draw[black,very thick] (3,5)--(13/3,14/3);
             \draw[black,very thick] (5,5)--(19/3,14/3);
               
                \draw[black,very thick] (6,2)--(7,7/4);
                 \draw[black,very thick] (6,4)--(7,15/4);
                 
                  \draw[very thick, fill=black]  (0,0) circle[radius=2.2pt];
                   \draw[very thick, fill=black]  (2,0) circle[radius=2.2pt];
                    \draw[very thick, fill=black]  (4,0) circle[radius=2.2pt];
                     \draw[very thick, fill=black]  (6,0) circle[radius=2.2pt];
                      \draw[very thick, fill=black]  (1,1) circle[radius=2.2pt];
                       \draw[very thick, fill=black]  (3,1) circle[radius=2.2pt];
                        \draw[very thick, fill=black]  (5,1) circle[radius=2.2pt];
                         \draw[very thick, fill=black]  (0,2) circle[radius=2.2pt];
                          \draw[very thick, fill=black]  (2,2) circle[radius=2.2pt];
                           \draw[very thick, fill=black]  (4,2) circle[radius=2.2pt];
                            \draw[very thick, fill=black]  (6,2) circle[radius=2.2pt];
                             \draw[very thick, fill=black]  (0,4) circle[radius=2.2pt];
                   \draw[very thick, fill=black]  (2,4) circle[radius=2.2pt];
                    \draw[very thick, fill=black]  (4,4) circle[radius=2.2pt];
                     \draw[very thick, fill=black]  (6,4) circle[radius=2.2pt];
                      \draw[very thick, fill=black]  (1,3) circle[radius=2.2pt];
                       \draw[very thick, fill=black]  (3,3) circle[radius=2.2pt];
                        \draw[very thick, fill=black]  (5,3) circle[radius=2.2pt];
                         \draw[very thick, fill=black]  (1,5) circle[radius=2.2pt];
                       \draw[very thick, fill=black]  (3,5) circle[radius=2.2pt];
                        \draw[very thick, fill=black]  (5,5) circle[radius=2.2pt];
                        
                         \draw[very thick, fill=white]  (1/3,2/3) circle[radius=2.5pt];
                         \draw[very thick, fill=white]  (7/3,2/3) circle[radius=2.5pt];
                         \draw[very thick, fill=white]  (13/3,2/3) circle[radius=2.5pt];
                         \draw[very thick, fill=white]  (19/3,2/3) circle[radius=2.5pt];
                         
                          \draw[very thick, fill=white]  (1/3,8/3) circle[radius=2.5pt];
                         \draw[very thick, fill=white]  (7/3,8/3) circle[radius=2.5pt];
                         \draw[very thick, fill=white]  (13/3,8/3) circle[radius=2.5pt];
                         \draw[very thick, fill=white]  (19/3,8/3) circle[radius=2.5pt];
                         
                          \draw[very thick, fill=white]  (1/3,14/3) circle[radius=2.5pt];
                         \draw[very thick, fill=white]  (7/3,14/3) circle[radius=2.5pt];
                         \draw[very thick, fill=white]  (13/3,14/3) circle[radius=2.5pt];
                         \draw[very thick, fill=white]  (19/3,14/3) circle[radius=2.5pt];
                         
                          \draw[very thick, fill=white]  (4/3,5/3) circle[radius=2.5pt];
                         \draw[very thick, fill=white]  (10/3,5/3) circle[radius=2.5pt];
                         \draw[very thick, fill=white]  (16/3,5/3) circle[radius=2.5pt];
                          \draw[very thick, fill=white]  (4/3,11/3) circle[radius=2.5pt];
                         \draw[very thick, fill=white]  (10/3,11/3) circle[radius=2.5pt];
                         \draw[very thick, fill=white]  (16/3,11/3) circle[radius=2.5pt];

      \end{tikzpicture}
      \caption{The image of the harmonic realization $\Phi_h$ in Example 2.}
     
      \end{figure}

\noindent Consider the zero range process on the discrete torus $ \mathbb{T}_N^2: =\mathbb{Z}^2 / N \Gamma$ with generator

\begin{equation}
L_N f ( \eta ):= \sum_{ x, y \in \mathbb{T}_N^2} p(x,y) g( \eta_x) [ f( \eta^{x,y})-f( \eta )].
\end{equation}
Note that weight function $p(\cdot, \cdot )$ is inhomogeneous, namely not invariant under the group action of $\mathbb{Z}^2$. As mentioned in Remark \ref{graph}, to investigate the hydrodynamic limit on $\mathbb{T}_N^2$, we can regard this process as a homogeneous process on the hexagonal lattice with the realization $\Phi$ in Ex 3b. More precisely, consider the zero range process on the hexagonal lattice, which is introduced in Section 2 Ex 3a, with generator 

\begin{equation}
L_Nf(\eta):=\sum_{e \in E_N} p(e) g(\eta_{oe}) [f(\eta^e) -f(\eta)]
\end{equation}
where the symmetric periodic weight function $p(\cdot)$ is given by 
\begin{equation}
p(e_1)=\frac{1}{3},\ \ p( e_2 )=\frac{ 1}{ 2 }, \ \ p( e_3 )= \frac{1 }{ 6}.
\end{equation} 
Here $e_1,e_2,e_3$ are the edges of the fundamental graph defined in Section 2 Ex 3a. Choosing the realization $\Phi$ in Ex 3b, the empirical process associated with $\Phi$ is  
 \begin{equation*}
  \pi_t^{\Phi,N} ( du )=\frac{1}{|V_N|} \sum_{x \in V_N}  \eta_x(t) \delta_{ \Phi_N (x)} (du).
 \end{equation*}
 Note that $\Phi$ is not harmonic for $p(\cdot)$ and the associated harmonic realization $\Phi_h$ is obtained by shifting every white vertex along (1/3,-1/3)(See Figure \ref{2}). By Theorem \ref{main}, we obtain that $\{ \pi_t^{\Phi ,N} \}_N $ converges to $  \rho(t, u) du/vol( \mathbb {T }_\phi^d ) $ in probability. Here $vol( \mathbb {T }_\phi^d ) =2$ and $\rho$ is the unique weak solution of $\partial_t\rho=\nabla \mathbb{D}_{\Phi_h} \nabla  \Psi(\rho)$, where $\mathbb{D}_{\Phi_h}$ is the diffusion matrix of $\Phi_h$,

 \begin{gather*} 
 \mathbb{D}_{\Phi_h}  = 
       \left (  
           \begin{array}{cc}   5/ 9  &  1/ 9 \\ 1/ 9 &  2 / 9            \end{array}  
       \right).
\end{gather*}

      \begin{remark}
      Both the above two examples are non-gradient systems, and if we apply the non-gradient method, the diffusion coefficients are not so clear since they are given by a variational formula. However, we can compute the diffusion coefficients explicitly by applying our main theorems.
      \end{remark}

 \begin{appendices}

      \section{Some lemmas on crystal lattices}

For a given lattice group

$$\phi(\Gamma)=\{ \sum_{i=1}^d k_i u_i\ | \ k_i\ integers\},$$

\noindent where $\{u_1,\cdots,u_d\}$ is a basis in $\mathbb{R}^d$. For $\textbf{x} \in \mathbb{R}^d$, \textbf{x} can be written uniquely as $\textbf{x}=\sum_{i=1}^d x_i u_i$. Define a norm on $\mathbb{R}^d$ by 

\begin{equation*}
|| \textbf{x}||_1:=\sum_{i=1}^d |x_i|,
\end{equation*}
for $\textbf{x}=\sum_{i=1}^d x_i u_i \in \mathbb{R}^d$. Define the distance $d_1$ by setting $d_1(\textbf{x},\textbf{x}^\prime)=||\textbf{x}-\textbf{x}^\prime||_1$ for $\textbf{x}, \textbf{x}^\prime \in \mathbb{R}^d$ and the induced metric in $\mathbb{T}^d_\phi$ from $d_1$ is also denoted by $d_1$.

Since $\Phi$ is periodic, there exists a constant $C_0:=\max_{x \in V} || \Phi([x]x_0)-\Phi(x)||_1<\infty$ such that

\begin{equation} \label{A1}
|| \Phi([x]x_0)-\Phi([z]x_0)||_1-2C_0 \le || \Phi(x)-\Phi(z)||_1 \le || \Phi([x]x_0)-\Phi([z]x_0)||_1+2C_0,
\end{equation}
for $x,z \in V$. Furthermore, since $|| \Phi([x]x_0)-\Phi([z]x_0)||_1=|[x]-[z]|$, (\ref{A1}) can be written as

\begin{equation}
| [x]-[z] |-2C_0 \le || \Phi(x)-\Phi(z)||_1 \le | [x]-[z] |+2C_0.
\end{equation}

Let $D_{\phi}:= \{ \sum_{i=1}^d t_i u_i\ |\ 0\le t_i<1,i=1,\dots,d \}$ be the fundamental parallelotope. For $\textbf{x} \in \mathbb{R}^d$, there exists a unique $\sigma_\textbf{x} \in \Gamma $ such that $\textbf{x} \in \phi(\sigma_\textbf{x})+D_\phi $. Define the map $[\cdot]: \mathbb{R}^d \to \Gamma$ by setting $[\textbf{x}]=\sigma_{\textbf{x}}$. Since $\Phi$ is periodic, there exists a constant $C_1:=\max_{\textbf{x} \in \mathbb{R}^d} || \Phi([\textbf{x}]x_0)-\textbf{x}||_1<\infty$ such that

\begin{equation}
|| \Phi([\textbf{x}]x_0)-\Phi([\textbf{z}]x_0)||_1-2C_1 \le || \Phi(\textbf{x})-\Phi(\textbf{z})||_1 \le || \Phi([\textbf{x}]x_0)-\Phi([\textbf{z}]x_0)||_1+2C_1,
\end{equation}
for $\textbf{x}, \textbf{z} \in \mathbb{R}^d$. Furthermore, we have that

\begin{equation}
| [\textbf{x}]-[\textbf{z}]|-2C_1 \le || \Phi(\textbf{x})-\Phi(\textbf{z})||_1 \le | [\textbf{x}] - [\textbf{z}] |+2C_1.
\end{equation}

For $\epsilon >0$, define the $\epsilon$-ball in $\mathbb{T}^d_\phi$ centered on $\textbf{z} \in \mathbb{T}^d_\phi$ by setting

\begin{equation*}
B_\textbf{z}(\epsilon):=\{ \textbf{x} \in \mathbb{T}^d_\phi\ |\ d_1(\textbf{x}, \textbf{z}) \le \epsilon\}.
\end{equation*}

\noindent and let $\chi_{\textbf{z},\epsilon}: \mathbb{T}^d_\phi \to \mathbb{R}$ be a characteristic function defined by

\begin{equation*}
\chi_{\textbf{z},\epsilon}:=\frac{vol(\mathbb{T}^d_\phi)}{vol(B_\textbf{z}(\epsilon))} 1_{B_\textbf{z}(\epsilon)}
\end{equation*}

\begin{lemma} \label{vol}
There exists a constant $C_3(\epsilon)>0$ depending only on $\epsilon$ such that for any $\textbf{z} \in \mathbb{T}^d_\phi$ and any $N \ge 1$, 
\begin{equation}
\left | \frac{vol(B_\textbf{z}(\epsilon))}{ vol(\mathbb{T}^d_\phi)} -\frac{ |\cup_{|\sigma| \le \epsilon N} \sigma D_{x_0} |} { |V_N|} \right| \le \frac{C_2(\epsilon)}{N},
\end{equation}
where $vol(A)$ stands for the volume of a Borel set $A$ and $|B|$ the cardinality of a set $B$. 
\end{lemma}

\noindent Proof: For any $\textbf{z} \in \mathbb{T}^d_\phi$, take a lift $\widetilde{\textbf{z}} \in \mathbb{R}^d$. For sufficiently small $\epsilon >0$, take a lift $B_{\widetilde{\textbf{z}}} (\epsilon) \subset \mathbb{R}^d$ of $B_\textbf{z} (\epsilon) \subset \mathbb{T}^d_\phi$. By (6.3) and (6.4), we have that 

\begin{equation*}
\bigcup_{ | \sigma | \le \epsilon N-2C_1} \sigma [\widetilde{\textbf{z}}] D_\phi \subset B_{\widetilde{\textbf{z}}} (\epsilon N) \subset \bigcup_{ | \sigma  | \le \epsilon N+2C_1} \sigma [\widetilde{\textbf{z}}] D_\phi,
\end{equation*}
which implies that

\begin{equation*}
\left |  N^d vol(B_{\widetilde{\textbf{z}}} (\epsilon )) -vol(\cup_{ |\sigma |\le \epsilon N} \sigma[\widetilde{\textbf{z}}] D_\phi )\right | \le vol(D_\phi) 2^d [ (\epsilon N+2C_1)^d-(\epsilon N-2C_1)^d].
\end{equation*}      

\noindent Note that $vol(B_{\widetilde{\textbf{z}}} (\epsilon ))=vol(B_{{\textbf{z}}} (\epsilon ))$, $vol(D_\phi)=vol(\mathbb{T}^d_\phi)$ and $|V_N|=N^d | D_{x_0}| $, we obtain that there exists a constant $C_2(\epsilon)$ depending only on $\epsilon$ such that

      \begin{equation*}
\left | \frac{vol(B_\textbf{z}(\epsilon))}{ vol(\mathbb{T}^d_\phi)} -\frac{ |\cup_{|\sigma| \le \epsilon N} \sigma D_{x_0} |} { |V_N|} \right| \le \frac{C_2(\epsilon)}{N}.
\end{equation*}
For the empirical density $\pi^N:=\frac{1}{|V_N|} \sum_{x \in V_N} \eta_x \delta_{\Phi_N (x)}$ on $\mathbb{T}^d_\phi$, $\eta \in Z_N$, we have the following lemma.
\begin{lemma} \label{chi}
There exists a constant $C_4(\epsilon)>0$ depending only on $\epsilon$ such that
\begin{equation}
\frac{1}{|\Gamma_N|} \sum_{\sigma \in \Gamma_N} \left | \langle \pi^N , \chi_{\Phi_N (\sigma x_0), \epsilon} \rangle -\bar{\eta}_{\sigma x_0, \epsilon N} \right | \le \frac{C_3(\epsilon)}{N} \frac{1}{|V_N|} \sum_{x \in V_N} \eta_x.
\end{equation} 
\end{lemma}

  \noindent Proof: For $\sigma \in \Gamma_N$, take a lift $\widetilde{\sigma x_0} \in V$ of $\sigma x_0$ and a lift $B_{(1/N)\Phi(\widetilde{\sigma x_0}) } (\epsilon) \subset \mathbb{R}^d$ of $B_{\Phi_N(\sigma x_0)}(\epsilon) \subset \mathbb{T}^d_\phi$. Similar to Lemma \ref{vol}, it holds that
  
  \begin{equation*}
  \bigcup_{ | \sigma^\prime |\le \epsilon N-2C_0 } \sigma^\prime \sigma D_{x_0} \subset \left\{x \in V\  \bigg|\  \left | \left | \frac{1}{N}\Phi(x)-\frac{1}{N} \Phi(\widetilde{\sigma x_0})  \right | \right |_1 \le \epsilon \right\} \subset  \bigcup_{ | \sigma^\prime |\le \epsilon N+2C_0 } \sigma^\prime \sigma D_{x_0}.
  \end{equation*}
Furthermore, take a lift $\widetilde{\eta} \in Z$ of $\eta \in Z_N$, then it holds that 

\begin{equation*}
\sum_{x \in V_N \atop \Phi_N(x) \in B_{\Phi_N(\sigma x_0)}(\epsilon)} \eta_x =\sum_{x \in V \atop \frac{1}{N} \Phi(x) \in B_{(1/N)\Phi(\widetilde{\sigma x_0}) } (\epsilon)} \widetilde{\eta}_x.
\end{equation*}      
      
     \noindent  We have that,
      
      \begin{equation*}
      \left| \frac{1}{|V_N|} \sum_{x \in V \atop || \frac{1}{N} \Phi(x)-\frac{1}{N}\Phi(\widetilde{\sigma x_0}) ||_1\le \epsilon } \widetilde{\eta}_x - \frac{1}{|V_N| } \sum_{x \in \sigma^\prime \sigma D_{x_0} \atop |\sigma^\prime| \le \epsilon N} \widetilde{\eta}_x \right | \le \frac{1}{|V_N|} \sum_{x \in \sigma^\prime \sigma D_{x_0} \atop \epsilon N-2C_0 \le |\sigma^\prime| \le \epsilon N +2C_0} \widetilde{\eta}_x.
      \end{equation*}
      
    \noindent  By Lemma \ref{vol}, it holds that
      
      \begin{align*}
      \left | \frac{vol(\mathbb{T}^d_\phi)} {|V_N | vol(B_{\Phi_N(\sigma x_0)}(\epsilon) ) } -\frac{1}{|\cup_{|\sigma^\prime | \le \epsilon N} \sigma^\prime \sigma D_{x_0} |} \right |  &\le \frac{vol(\mathbb{T}^d_\phi)}{ vol( B_{\Phi_N(\sigma x_0)}(\epsilon)) |\cup_{|\sigma^\prime | \le \epsilon N} \sigma^\prime \sigma D_{x_0} |} \frac{C_3(\epsilon)}{ N} \\
     & \le \frac{C_5(\epsilon)}{N^{d+1}},
      \end{align*}
   where $C_5(\epsilon)$ is a constant depending only on $\epsilon$. By triangular inequality, we have that
   \begin{align*}
   \left | \langle \pi^N , \chi_{\Phi(\sigma x_0), \epsilon} \rangle -\bar{\eta}_{\sigma x_0, \epsilon N} \right | &\le \frac{vol(\mathbb{T}^d_\phi)} {|V_N | vol(B_{\Phi_N(\sigma x_0)}(\epsilon) ) } \sum_{x \in \sigma^\prime \sigma D_{x_0} \atop \epsilon N-2C_0 \le |\sigma^\prime| \le \epsilon N +2C_0} \eta_x 
   \\
   &\ \ +\left | \frac{vol(\mathbb{T}^d_\phi)} {|V_N | vol(B_{\Phi_N(\sigma x_0)}(\epsilon) ) } \sum_{x \in \sigma^\prime \sigma D_{x_0} \atop |\sigma^\prime| \le \epsilon N} \eta_x - \frac{1}{|\cup_{|\sigma^\prime | \le \epsilon N} \sigma^\prime \sigma D_{x_0} |} \sum_{x \in \sigma^\prime \sigma D_{x_0} \atop |\sigma^\prime| \le \epsilon N} \eta_x\right|
   \\
   &\le  \frac{vol(\mathbb{T}^d_\phi)} {|V_N | vol(B_{\Phi_N(\sigma x_0)}(\epsilon) ) } \sum_{x \in \sigma^\prime \sigma D_{x_0} \atop \epsilon N-2C_0 \le |\sigma^\prime| \le \epsilon N +2C_0} \eta_x +\frac{C_5(\epsilon)}{N^{d+1}} \sum_{x \in \sigma^\prime \sigma D_{x_0} \atop |\sigma^\prime| \le \epsilon N} \eta_x.
   \end{align*}    
   
   \noindent Sum on $\sigma \in \Gamma_N$, it concludes that there exists a constant $C_4(\epsilon)$ depending only on $\epsilon$ such that 
   
   \begin{equation*}
\frac{1}{|\Gamma_N|} \sum_{\sigma \in \Gamma_N} \left | \langle \pi^N , \chi_{\Phi(\sigma x_0), \epsilon} \rangle -\bar{\eta}_{\sigma x_0, \epsilon N} \right | \le \frac{C_3(\epsilon)}{N} \frac{1}{|V_N|} \sum_{x \in V_N} \eta_x.
\end{equation*}

 \end{appendices}

\bibliography{references}

\end{document}